\documentclass[USenglish]{article} 
\usepackage[small]{dgruyter}
\usepackage[utf8]{inputenc}
\usepackage{amsmath, amsthm, amscd, amsfonts,amssymb,graphicx,enumerate}
\numberwithin{equation}{section}
\usepackage{color}

\makeatletter
\@namedef{subjclassname@2010}{%
  \textup{2010} Mathematics Subject Classification}
\makeatother

\def\pic #1 by #2 (#3){\vbox to #2{\hrule width 
#1 height 0pt depth 0pt\vfill\special{picture #3}}}
\def\scaledpicture#1
by #2 (#3 scaled #4){{\dimen0=#1 \dimen1=#2\divide\dimen0 by
 1000 \multiply\dimen0 by #4\divide\dimen1 by 1000 \multiply\dimen1 
by #4\pic \dimen0 by \dimen1 (#3 scaled #4)}}

\newcommand{\comm}[1]{{}}
\newcommand{\sB}{{\tsp {B}}}

\renewcommand{\phi}{\varphi}
\newcommand{\Fr}{\text{\rm Fr}}

\newcommand{\NRC}{\text{\bf NRC}}
\newcommand{\PAC}{\text{\bf PAC}}

\newcommand{\C}{{\text{\rm C}}}

\newcommand{\by}{{\pmb y}}
\newcommand{\tP}{{\text{\rm P}}}

\begin{document}
\baselineskip=17pt
\hoffset-.25in



\newcommand{\tsec}{\arabic{section}}
\newcommand{\tsubsec}{\tsec.\arabic{subsection}}
\newcommand{\tsubsubsec}{\tsubsec.\arabic{subsubsection}}
\newcommand{\tsubsubsubsec}{\tsubsubsec.\roman{subsubsubsection}}
\newcommand{\teq}{\arabic{section}.\arabic{equation}}
\newcommand{\ttopic}{\arabic{section}.\arabic{topic}}
\newcommand{\teql}{\Alph{section}.\arabic{equation}}
\newcommand{\unfinished}{} 




\renewcommand{\theenumi}{\alph{enumi}}
\newcommand{\thelabelenumi}{\teq\alph{enumi}}

\newcommand{\sqr}[2]{{\vcenter{\vbox{\hrule height.#2pt\hbox{\vrule width.#2pt
height#1pt \kern#1pt\vrule width.#2pt}\hrule height.#2pt}}}}
\newcommand\msquare{\mathchoice\sqr34\sqr34\sqr{2.1}3\sqr{1.5}3}
\newcommand{\ssquare}{{\qquad\hfill$\square$}}

\newenvironment{explanation}{\noindent \bf Explanation:\rm\ }{\ssquare}
\newcounter{eqcount}
\newcounter{ttopic}
\renewcommand{\theeqcount}{\alph{eqcount}}
\renewcommand{\labelenumi}{{{\rm (\teq \alph{enumi})}}} 
\newenvironment{edesc}{\refstepcounter{equation}\begin{enumerate}}%
{\end{enumerate}}
\newenvironment{topic}{\refstepcounter{ttopic}\begin{topic}%
{{\hbox{\rm(\teq)}}} \item }{\end{topic}} 
\newenvironment{trivl}{\refstepcounter{equation}\begin{list}%
{{\hbox{\rm(\teql)\ }}} \\item }{\end{list}}
\newcommand{\display}[1]{{\vskip.1in\noindent {\eu  #1:}}} 
\newcommand{\sdisplay}[1]{{{$\scriptstyle\bullet$}\! #1 \!{$\scriptstyle\bullet$}}} 
\newcommand{\Sdisplay}[1]{{\vskip.2in\noindent $\star$ {\eu #1} $\star$\hskip.1in:}} 

\newcommand{\ring}[1]{{\mathbb #1}}
\newcommand\bZ{{\ring{Z}}}
\newcommand{\field}[1]{{\mathbb #1}}
\newcommand\bC{{\ring{C}}} \newcommand\bR{{\ring{R}}}
\newcommand\bF{{\ring{F}}} \newcommand\bQ{{\ring{Q}}}
\newcommand\bH{{\ring{H}}}
\newcommand{\csp}[1]{{\mathbb #1}}
\newcommand{\tsp}[1]{{\mathcal #1}}
\newcommand{\stsp}[1]{{\text{\sCal #1}}}
\newcommand{\esp}[1]{{\mathcal #1}}
\newcommand{\sesp}[1]{{\mathcal #1}}
\newcommand{\bM}{\csp{M}}
\newcommand{\bB}{\csp{B}}
\newcommand{\prP}{\csp{P}}
\newcommand{\afA}{\csp{A}}
\newcommand{\sA}{{\esp{A}}} \newcommand{\ssA}{{\sesp{A}}}
\newcommand{\sC}{{\tsp{C}}} \newcommand{\sU}{{\tsp{U}}}
\newcommand{\sO}{{\tsp{O}}} \newcommand{\sI}{{\tsp{I}}}
\newcommand{\sQ}{\tsp{Q}}
\newcommand{\sP}{{\tsp {P}}} \newcommand{\ssP}{{\stsp {P}}}
\newcommand{\sL}{{\tsp {L}}} \newcommand{\sS}{{\tsp {S}}}
\newcommand{\sT}{{\tsp {T}}} \newcommand{\sH}{{\tsp {H}}}
\newcommand{\sX}{{\tsp {X}}} \newcommand{\sY}{{\tsp {Y}}}
\newcommand{\sM}{{\tsp {M}}} \newcommand{\sD}{{\tsp {D}}}
\newcommand{\sG}{{\tsp {G}}} \newcommand{\sE}{{\tsp {E}}}
\newcommand{\sR}{{\tsp {R}}} \newcommand{\GT}{{$\widehat{\tsp {GT}}$}}
\newcommand{\bD}{{\csp {D}}}
\newcommand{\bT}{{\csp {T}}} \newcommand{\bN}{{\csp {N}}}

\newcommand{\La}{{\Lambda}}

\newcommand{\eql}[2]{{\rm (\ref{#1}\ref{#2})}} 

\newcommand{\vect}[1]{{\pmb #1}} \newcommand{\ori}{{\pmb 0}}
\newcommand{\ba}{{\vect{a}}} \newcommand{\bg}{\vect{g}}
\newcommand{\bd}{{\vect{d}}} \newcommand{\be}{{\vect{e}}}
\newcommand{\bp}{{\vect{p}}} \newcommand{\bx}{{\vect{x}}}
\newcommand{\bv}{{\vect{v}}} \newcommand{\bw}{{\vect{w}}}
\newcommand{\bs}{{\vect{s}}} \newcommand{\bz}{{\vect{z}}}
\newcommand{\bh}{{\vect{h}}}  \newcommand{\bu}{{\vect{u}}}
\newcommand{\row}[2]{{#1_1,\ldots,#1_{#2}}}
\newcommand{\rowb}[3]{{#1_{#2},\ldots,#1_{#3}}}
\newcommand{\smatrix}[4]{{\big(\begin{array}{cc}
\!\lower2pt\hbox{$\scriptstyle#1$} &\lower2pt\hbox{$\scriptstyle#2$}\!
\\\! \raise2pt\hbox{$\scriptstyle#3$} &\raise2pt\hbox{$\scriptstyle#4$}
\!\end{array}\big)}}
\newcommand{\col}[2]{{\big(\begin{array}{c}
\!\lower2pt\hbox{$\scriptstyle#1$}  \!
\\\! \raise2pt\hbox{$\scriptstyle#2$}
\!\end{array}\big)}}

\newcommand{\scases}[5]{{#1 =\left\{\begin{array}{ll} #2 &
\mbox{for $#3$}\\ #4 & \mbox{for $#5$}.\end{array}\right.}}
\newcommand{\Scases}[7]{{#1 =\left\{\begin{array}{ll} #2 &
\mbox{for $#3$}\\ #4 & \mbox{for $#5$}\\#6 & \mbox{for
$#7$}.\end{array}\right.}}
\newcommand{\SScases}[9]{{#1 =\left\{\begin{array}{ll} #2 &
\mbox{for $#3$}\\ #4 & \mbox{for $#5$}\\#6 & \mbox{for $#7$}
\\ #8 & \mbox{for $#9$}.\end{array}\right.}}

\newcommand{\texto}[1]{{\textr{#1}}}
\newcommand{\GL}{\texto{GL}} \newcommand{\SL}{\texto{SL}}
\newcommand{\SO}{\texto{SO}} \newcommand{\ind}{\texto{ind}}
\newcommand{\PSL}{\texto{PSL}} \newcommand{\PGL}{\texto{PGL}}
\newcommand{\Gal}{\texto{Gal}} \newcommand{\sni}{\texto{SNi}}
\newcommand{\Hom}{\texto{Hom}} \renewcommand{\ni}{\texto{Ni}}
\newcommand{\Spec}{\texto{Spec}} \newcommand{\Pic}{\texto{Pic}}
\newcommand{\textr}[1]{{\text{\rm #1}}}
\newcommand{\tr}{\textr{tr}} \newcommand{\ord}{\textr{ord}}
\newcommand{\abs}{\textr{abs}}  \newcommand{\cyc}{\textr{cyc}}
\newcommand{\proj}{\textr{proj}} \newcommand{\spn}{\textr{span}}
\newcommand{\pC}{{\textr{C}}} \newcommand{\inn}{\textr{in}}
\newcommand{\diag}[1]{{\textr(#1)}} \newcommand{\Aut}{\textr{Aut}}
\newcommand{\pr}{\textr{pr}}
\newcommand{\ci}{{i\,}}

\newcommand{\norm}{{\triangleleft\,}}
\newcommand{\RET}{{\text{\rm RET}}}
\newcommand{\BCL}{{\text{\rm BCL}}}
\newcommand{\IGP}{{\text{\rm IGP}}}
\newcommand{\ext}{\texto{ext}}
\newcommand{\alg}{\texto{alg}}
\newcommand{\abe}{\texto{ab}}
\newcommand{\rd}{\texto{rd}}
\newcommand{\ari}{\texto{ar}}
\newcommand{\nil}{\texto{nil}}
\newcommand{\tG}[1]{{}_{#1}\tilde G}
\newcommand{\Syp}{\text{\rm Sp}}
\newcommand{\Det}{\text{\rm Det}}
\newcommand{\bbQ}{{\bar{\ring{Q}}}}
\newcommand{\GAP}{{\bf GAP}}
\newcommand{\ot} {\otimes}

\newcommand{\sph}{{\vphantom 1}}

\newcommand{\textb}[1]{{\text{\bf #1}}}
\newcommand{\bfC}{{\textb{C}}}
\newcommand{\longmapright}[2]{\smash{\mathop{\hbox to
#2pt{\rightarrowfill}}\limits^{#1}}}
\newcommand{\Longmapright}[2]{\smash{\mathop{\hbox to
#2pt{\Rightarrowfill}}\limits^{#1}}}
\newcommand{\longmapleft}[2]{\smash{\mathop{\hbox to
#2pt{\leftarrowfill}}\limits^{#1}}}
\newcommand{\mapdown}[1]{\Big\downarrow\rlap{$\vcenter{\hbox{$\scriptstyle#1$}}
$}} \newcommand{\lmapdown}[1]{\llap{$\vcenter{\hbox{$\scriptstyle{#1}$}}
$}\Big\downarrow}
\newcommand{\mapup}[1]{\Big\uparrow\rlap{$\vcenter{\hbox{$\scriptstyle#1$}}
$}} \newcommand{\lmapup}[1]{\llap{$\vcenter{\hbox{$\scriptstyle{#1}$}}
$}\Big\uparrow}
\newcommand{\mapright}[1]{\smash{\mathop{\longrightarrow}\limits^{#1}}}
\newcommand{\mapleft}[1]{\smash{\mathop{\longleftarrow}\limits^{#1}}}
\newcommand{\st}[2]{\stackrel{#1}{#2}}

\newcommand{\np}{{+}}   \newcommand{\nm}{{-}}
\newcommand\sem{\setminus}
\newcommand{\dd}[2]{{\frac{d#1}{d#2}}}
\newcommand{\lrang}[1]{{\langle #1\rangle}}
\newcommand{\blrang}[1]{{\big< #1\big>}}
\newcommand{\con}[3]{{#1 \equiv #2 \bmod #3}}
\newcommand{\Div}{{\,|\,}}
\newcommand{\eqdef}{\stackrel{\text{\rm def}}{=}}
\newcommand{\bsl}{\backslash}
\newcommand{\pa}[2]{{\frac{\partial #1} {\partial #2}}}


\newcommand{\fig}[3]{{ \begin{figure}\centering
\includegraphics{#1}\caption{#2} \label{#3}\end{figure}}}

\newcommand{\nvect}[3]{\begin{picture}(0,0)\vector(#1,#2){#3}\end{picture}}
\newcommand{\hoval}[3]{\begin{picture}(0,0)\oval(#1,#2)[#3]\end{picture}}
\newcommand{\qoval}[4]{\begin{picture}(0,0)\oval(#1,#2)[#3#4]\end{picture}}

\newfont{\sevenrm}{cmr7}
\newfont{\bsevenrm}{cmbx7}
\newfont{\mathseven}{cmsy7}
\newfont{\bigmath}{cmsy10 scaled 1200}
\newfont{\fiverm}{cmr5}
\newfont{\bfiverm}{cmbx5}
\newfont{\hel}{cmbx10 scaled 1400}
\newfont{\eu}{eufb10}
\newfont{\sseu}{eufm5}
\newfont{\seu}{eufm7}
\newfont{\Cal}{cmmib10}
\newfont{\sCal}{cmmib7}
\newfont{\zch}{eusb10}


\newtheorem{theorem}{Theorem}[section]

\theoremstyle{plain}
\newtheorem{thm}{Theorem}[section] 
\newtheorem{lem}[thm]{Lemma}
\newtheorem{princ}[thm]{Principle}
\newtheorem{fact}[thm]{Fact}
\newtheorem{prop}[thm]{Proposition}
\newtheorem{cor}[thm]{Corollary}
\newtheorem{goal}[thm]{Goal}
\newtheorem{res}[thm]{Result}
\newtheorem{thml}{Theorem} %
\newtheorem{leml}[thml]{Lemma}
\newtheorem{princl}[thml]{Principle}
\newtheorem{factl}[thml]{Fact}
\newtheorem{propl}[thml]{Proposition}
\newtheorem{corl}[thml]{Corollary}
\newtheorem{goall}[thml]{Goal}
\newtheorem{resl}[thml]{Result}


\theoremstyle{definition}
\newtheorem{defn}[thm]{Definition}
\newtheorem{exmp}[thm]{Example}
\newtheorem{guess}[thm]{Conjecture}
\newtheorem{quest}[thm]{Question}
\newtheorem{prob}[thm]{Problem}
\newtheorem{questl}[thml]{Question}
\newtheorem{defnl}[thml]{Definition}
\newtheorem{exte}[thm]{Ext.}
\theoremstyle{remark}
\newtheorem{rem}[thm]{Remark}
\newtheorem*{sol}{Solution}

\newcommand{\xs}{\times^s\!}
\newcommand{\wsp}{{$\,$---$\,$}} 
\newcommand{\lp}{{\rm(}} \newcommand{\rp}{{\rm)}}


\newcommand{\twid}{\~{\hbox{$\,$\!\!}}}

\newcommand{\data}[1]{\hskip3truein\hbox{\vtop{\hbox{Michael Fried}
\hbox{Math Dept: UC Irvine}
\hbox{Home Phone:  714-854-3634}
\hbox{#1}
\hbox{\sevenrm e-mail: mfried@math.uci.edu}}}}

\newenvironment{exmpl}{\begin{exmp}}{\hfill $\triangle$ \end{exmp}}

\title{\hskip .25 in Diophantine statements  over Residue fields:  \\ \hskip .5in Galois Stratification and Uniformity } 


\date{} 
\maketitle 

\begin{center} Michael D.~Fried, Emeritus Professor of Mathematics \\ University of California at Irvine\end{center} 

\begin{abstract}\ We consider three general problems about Diophantine statements   over finite fields that connect to the {\sl Galois stratification\/}  procedure for deciding such problems. To bring to earth the generality of these problems, we clarify each using the negative solution of U.~Felgner's simply stated question relating pairs of finite fields $\bF_p$ and $\bF_{p^2}$ for large primes $p$.  

For a diophantine problem, $D$, interpretable for almost all primes $p$, the paper plays on attaching a Poincar\'e series, $P_{D,p}$, to $(D,p)$. The rationality and computability of $P_{D,p}$ as $p$ varies gives some aspects of continuity. Most interesting, though, is how the coefficients of $P_{D,p}$   vary with $p$. The paper notes three ways to give counts for those coefficients:  

\begin{itemize} \item[$\bullet$] points mod $p$ achieving particular conjugacy classes in a Galois stratification; 
\item[$\bullet$] points mod $p$ on absolutely irreducible varieties; or 
\item[$\bullet$] traces of the $p$ Frobenius on a Chow motive. \end{itemize} 
\hskip .25 in The gist of this paper, using one explicit well-known diophantine problem,  is that  Galois stratification  naturally tethers these three abstract approaches.   Among \cite{FrS76}, \cite{FrHJ84}, \cite{DL98}, and \cite{FrJ04} is a history of the fundamental ideas and original motivations for extending the Galois Stratification procedure. Others, \cite{W03}, \cite{Hr12} and \cite{To16}  referred to in \S\ref{motivations} show a more extensive set of motivations for extending the topic today. So, it is time for an introduction relating the ideas to Deligne's proof of the Weil Conjectures and Langland's Program applied to two specific recognizable problems/examples. 
 
\end{abstract} 

\section{Three general problems on finite field equations} Traditional diophantine statements consider algebraic subsets of affine space with blocks of variables with some, say  $\bx=(\row x m)$, unquantified, and others, say $\by=(\row y n)$, quantified by the symbols $\exists$ and $\forall$. Questions over finite fields often assume the equation coefficients are in $\bZ$. For example, they might consider  $\bx$ and $\by$ taking values in the residue class fields $\bF_p$, and ask if the statements are true for almost all $p$. 

\subsection{Introducing the three problems}  \label{introduction} 
Long time results on variants using {\sl Galois stratification\/} (\S\ref{galstrat}) have detailed literature as in \cite[Chaps.~30-31]{FrJ04}.  There are, however, three topics inadequately treated there.  
\begin{edesc} \label{quest} \item \label{questa} Testing how a quantified $D$ over distinct residue fields cohere -- {\sl uniformity in $p$\/} -- using equivalent data from reducing an unquantified object over $\bQ$  mod most primes $p$.  
\item  \label{questb} Coordinating (*) affine space arguments and (**) arithmetic homotopy with scheme or {\sl projective geometry\/} language (the natural domain of arithmetic geometry). 
\item \label{questc} Considering statements with variables taking values in the algebraic closure of $\bF_p$, but fixed by respective powers of the Frobenius $\Fr_p$: \begin{equation} \text{for }\bx \text{ (resp.~$\by$)}(\Fr_p^{d_1},\dots, \Fr_p^{d_m}) \eqdef \Fr_p^{\bd} \text{ (resp.~} (\Fr_p^{e_1},\dots, \Fr_p^{e_n})=\Fr_p^{\be}). \end{equation}   
\end{edesc} 
\hskip .25 in To help understand  aspects of all three problems we use Felgner's Problem:  $$\text{ Can we define the fields $\bF_p$ within the theory of the  fields $\bF_{p^2}$. }$$

\display{Comments on notation and background} We assume the reader knows the meaning of the words {\sl cover of normal quasi-projective varieties\/}, where cover may, usually here, includes ramification, but is a {\sl finite, flat\/} morphism (Grothendieck's definition). Although \cite{H77} is much more comprehensive, we still find it valuable to refer to \cite{Mu66} for proofs through special illuminating cases. Especially for Segre's proof that normalization in a function field of a quasi-projective variety is also quasiprojective \cite[p.~400]{Mu66}. The constructions alluded to in this paper also depend on Chevalley's Theorem that {\sl the image of a constructible set is constructible\/} \cite[p.~97]{Mu66}. 

Neither is illuminating on the phrase {\sl non-regular Galois cover\/} for a very good reason; Both do almost everything over algebraically closed fields. The Galois stratification procedure confronts situations in which $\phi: Y\to X$ is a cover of absolutely irreducible varieties over a (perfect) field $K$ (with algebraic closure $\bar K$), while its Galois closure $\hat \phi: \hat Y \to X$, $\hat Y$ is not absolutely irreducible over $K$. In other words, $K(\hat Y)\cap \bar K$ is a nontrivial extension of $K$.

\display{Comments on \eql{quest}{questa}} For any prime $p$ in \eql{quest}{questa}, a {\sl Galois stratification\/} allows eliminating the quantifers on  the variables, producing a quantifer free statement. The $\bQ$ object in \eql{quest}{questa} -- a Galois stratification object, but over $\bQ$ rather than over a finite field -- has not been previously emphasized, though it has always been present. Though more general than a pure algebraic variety, it allows all the elimination theory and reductions mod primes. 

At the heart of the Galois stratification procedure is its use of the {\sl non-regular Chebotarev analog\/}. We use the acronym \NRC. Both the non-regular adjective and the Chebotarev error term impact its use, as in the negative solution of Felgner's Problem. It appears in two ways: 
\!\!\!\begin{edesc} \label{pieces} \item \label{piecesa} Eliminating quantifiers, the main reason we introduced the stratification procedure; and  
\item \label{piecesb} estimating the unquantified variable values satisfying the diophantine conditions. 
\end{edesc} 
\hskip .25 in We can untangle the two steps of \eqref{pieces}. So doing reveals two separate aspects of the stratification procedure. Its use in \eql{pieces}{piecesa} for collecting finite fields points. Then,  its generalization of quantifier elimination over many other collections of fields, as appears in \cite[Chap.~30]{FrJ04} as well as the variant given in \eql{quest}{questc}. 

Starting with an elementary diophantine statement (problem) $D$, the procedure generalizes elementary statements to Galois stratifications $\sS$. In each generalizing case, every step of the stratification procedure produces a new Galois stratification. The number of steps depends on the number of blocks of quantified variables.   \S\ref{galstrat} gives the notation and explains the main theorem as given in several sources starting from \cite{FrS76}. 

That is, suppose there are $k$ blocks of quantified variables. For each prime $p$ this produces $\bar \sS_p\eqdef\{\sS_{p,k}, \dots, \sS_{p,0}\}$, a sequence of Galois stratifications (starting at $k$ going to 0), and a finite Galois extension $\hat K/\bQ$ (with group $G_P$). Excluding a finite computable set of primes $p$ these have the following property. For $p$ a prime unramified in $\hat K$, denote an element of the  conjugacy class of the Frobenius at $p$ by $\Fr_p$, and its fixed field in $\hat K$ by $\hat K_{p}$ . 
  \begin{equation} \label{stratproc} \begin{array}{c} \text{For each conjugacy class $\C$ of values of the Frobenius, $\Fr_p\in G_P$,} \\ \text{there is an object $\bar \sS_\C$ over $\hat K_{p}$ whose reduction mod a prime above $p$ gives $\bar \sS_p$.}\end{array} \end{equation}   
That gives a start to many aspects of uniformity in $p$. 

One form of the stratification procedure moves through field theory -- after all, {\sl Field Arithmetic\/} is the \cite{FrJ04} title --  based on {\sl Frobenius fields}. While we still rely on details on the stratification procedure, \S\ref{quoteFrJ04} gives an overview hitting the most significant of those details. Especially it shows the role of the Chebotarev analog. 

\S\ref{nonregcheb} emphasizes the word {\sl non-regular\/} in that analog   and its subtleties.  For example, this gives the main distinctions, running over classes $\C$ in $G$, between the collections $\bar \sS_\C$.  
  
\display{Comments on \eql{quest}{questb}}  From the beginning there was an idea of attaching a more concise object, a Poincare series $P_{p}\eqdef P_{D,p}$, to  $\bar \sS_p$. The stratification procedure was the main tool in showing this object to be  a computable rational function. The tension between (*) and (**) in \eql{quest}{questb} appeared long before  \cite{FrS76} as recounted in \cite{Fr12} where this author's applications to finite field questions used Riemann's work and projective geometry. The tension continued: \begin{itemize} \item[$\bullet$] on one hand, in applying Bombieri's use of Dwork's affine/singular theory of the zeta function to explicitly compute stratification Poincar\'e series {\sl characteristics\/};  and, 
\item[$\bullet$]  on the other, the Denef-Loeser (and Nicaise) production of Poincar\'e series coefficients as  Chow motives. Deligne's proof of the Weil conjectures is in the background. \end{itemize} 

\hskip .25 in Continuing the comments on \eql{quest}{questa}, the {\sl nonregular\/} adjective shows in Thm.~\ref{nonregcover}. It gives how the stratification procedure varies with either the prime $p$ (or in applying it to a specific finite field, the extension $\bF_q$ for $q$ a power of $p$).

 \S\ref{concfelgner} applies the Main Thm. to Felgner's problem. Here -- I suspect a careful reader will agree -- a mere cardinality estimate, rather than some precise understanding of the variance with either the prime $p$ or the extension $\bF_q$ for $q$ a power of $p$, seems inadequate. 
 
 Yet, there is considerable flexibility in how those Poincar\'e series capture information.  \S\ref{schurcover} illustrates with an a diophantine example the author has used with undergraduates. (Proofs are another matter.) This example helps picture how the reference to Chow motives, as pieces of \'etale cohomology groups,  works. 

\display{Comments on \eql{quest}{questc}} \S\ref{Felgner} explains Felgner's problem in very down-to-earth language. Felgner  is a piece from the case $m=1$, $d_1=1$, with all  $e_i\,$s equal to 2.   It poses for Galois stratifications -- expanded by {\sl Frobenius vectors\/} (like $\Fr_q^{\bd}$)  -- the nature of corresponding Poincar\'e series. For example, when can they be derived from standard Poincar\'e series?  

This latter story passes through work of D.~Wan, who introduced a Zeta function in the unquantified case, generalizing Artin's zeta where all $d_i\,$s are 1, and Hrushovski and Tomasic who took a model theoretic approach, enhanced by Galois stratification.\footnote{This starts with Ax and Kochen on the Artin Conjecture on $p$-adic points on forms of degree $d$ in projective $d^2$-space. To understand my motivation, consider the down-to-earth problems I connected to classical arithmetic geometry as recounted in the partially expository \cite{Fr05} and \cite{Fr12}.}

From constraints on space and time, we limit our responses to questions \eql{quest}{questb}  and \eql{quest}{questc} to merely motivating seriously using these Poincar\'e series, as if their coefficients are Chow motives to which you apply the Frobenius $\Fr_p$. Our goal is to ask questions about their variance with $p$. Especially for what $p$ certain especially simple results occur, as given explicitly in our example.  We expect to  expand  on parts  of \S\ref{hrushwan} in a later paper. 

\subsection{Felgner clarifies what are statements over finite fields} \label{Felgner} 

A traditional question  in the theory of finite fields starts with blocks of variables $\by_i\in \afA^{n_i}$, $i=1,\dots,k$. Below we use the notation $\sum_{j=1}^k n_j\eqdef N_{k}$.  Assume these blocks are respectively quantified by $\row Q k={\pmb Q}$. Also we have a separate unquantified block of variables $\bx\in \afA^m$ together with an algebraic set $X$ described by  
$$\text{$\{(\bx, \row \by k)\mid \phi(\bx,\row {\by}k)\}$.}$$ 

Here $\phi$ is a collection of polynomial equalities and inequalities with coefficients in a ring $R$, defining a constructible subset of $\afA^{m\np N_k}$. A traditional question is, given $\phi$ and $\row Q k$, and no $\bx$ variables, is there a useful procedure for deciding if the statement is true for {\sl almost all\/} finite fields $\bF_q$ in the following collections.  
\begin{edesc} \label{tradff} \item  $R=\bF_p$: $q$ are powers of a given $p$. 
\item  $R=\sO_K$, ring of integers of a number field $K$: $q=|\sO_K/\bp|$, orders of residue class fields. 
\item  $R=\bZ$, the whole Universe case: All $\bF_q$. 
\end{edesc} 

\S\ref{galstrat} explains the Main Theorem of Galois stratification. It is  a stronger result  than that above in that our starting statement usually includes unquantified variables $\bx$. Then, the result is that the starting statement -- with, excluding $\bx$, all variables quantified  -- is equivalent to one with no quantified variables. That is, it gives an {\sl elimination of quantifiers}.  We use Felgner's problem \cite{Fe90} to show its value and how it works. 

In each case, in changing \lq\lq almost all\rq\rq\ to \lq\lq all,\rq\rq,  our method (and generalizations) leaves checking a finite number of prime powers. For some problems those exceptional $q\,$s could be accidents. For others -- motivated by classical cases -- they could be significant. The Poincar\'e series approach (\S\ref{poincare}), seeking a significant rational function, has been the gold standard.   

Felgner -- in language reminiscent of the above (with $m=1$) -- asks if given $x'\in \bF_{p^2}$, 
$$\text{is  $(x',\row {\by}k)\in X$  true for $Q_1\by_1\dots Q_k\by_k\in (\bF_{p^2})^{N_k}$,  if and only if $x'\in \bF_p$? }$$
Except we are asking about a different set -- $\{\bF_{p^2}\mid p  \text{ prime }\}$ -- of fields than those included in any of the lists \eqref{tradff}.  Our method easily  tolerates such flexibility. Further, the negative conclusion holds even if we asked for its truth replacing $p$ by prime-powers $q$ ($\bF_p\mapsto \bF_q$, $\bF_{p^2}\mapsto \bF_{q^2}$) and just for {\sl infinitely\/} many $q$. 

The conclusion about the set of such $x'$ opens territory that isn't hinted at by the list \eqref{tradff}. Still, the nature of that set arising with actual problems in finite fields (as in \S\ref{unifp}),  is the model for which the original result aimed.  
Here is an example of what we will exclude. Take $k=2$, and $\pmb Q=(\exists,\forall)$. Read this as follows. 
\begin{edesc}  \label{Felgnersimple}  
\item Given $x'\in \bF_{q}$, there exists $\by_1'\in (\bF_{q^2})^{n_1}$ with  $\by_2'\in (\bF_{q^2})^{n_2}$, so that $(x',\by_1',\by_2')\in X$. 
\item But if $x'\in \bF_{q^2}\setminus \bF_q$, then for any $\by_1'\in (\bF_{q^2})^{n_1}$  there is $\by_2'\in (\bF_{q^2})^{n_2}$ with $(x',\by_1',\by_2')\not \in X$. 
 \end{edesc}  

\begin{center}  It seems so simple: $x\in \bF_{q^2}$ is in $\bF_q$ if and only if $\Fr_q(x)\eqdef x^q=x$. \\ How could it fail that some algebraic statement generalizing \eqref{Felgnersimple} would encapsulate this? \end{center} 

\section{The Main Theorem of Galois Stratification} The main theorem inductively, \S\ref{galstrat}, eliminates quantifiers for problems that 
\begin{equation} \label{elemstat} \begin{array}{c} \text{start as elementary statements interpretable over {\sl any\/} extension}\\ \text{ of a residue class field of the ring of integers of a number field $K$.}\end{array} \end{equation}  
First we give the appropriate Main Theorem version that allows the elimination. This includes all possible statements  that might have answered Felgner's question in the affirmative. To simpify discussing the procedure, take as fundamental the following problem.
$$\text{Is a statement true for {\sl almost all\/} residue class fields $\bZ/p=\bF_p$ (that is, $R=\bZ$ in \eqref{tradff}).}$$ 

Given the elementary statement start of \eqref{elemstat}, we expediently consider the {\sl exceptional stratification\/} primes (Def.~\ref{excstratprimesdef}) not included in the phrase \lq\lq almost all.\rq\rq\footnote{The word {\sl exceptional\/} has two, incompatible uses in this paper, forced on us by the historical record. I have put the word {\sl stratification\/} in this one to separate them.} That thereby allows the procedure to consider if a statement holds for {\sl all\/} primes or variants of that question. 

The philosophical differences between \cite{FrS76} and \cite{FrJ04} are small, though there are more details in the latter. That  is primed to handle -- {\sl Frobenius field\/} --  results more general than the finite field case. It shows there is a {\sl general Galois stratification\/} idea. Then, you adjust by using an \NRC\ replacement for each collection of (so-called) Frobenius fields. 

This approach produces the finite field case by observing, if $q$ is large, then the collection of finite fields behaves like a collection of Frobenius fields.\footnote{That's no surprise since that was in the original motivations for \cite{FrJ04}.} We don't redo that.  Rather we show a path through  \cite[Chap.~30]{FrJ04} in  \S\ref{quoteFrJ04} -- with an example -- that works.  

\subsection{Galois stratification and elimination of quantifiers} \label{galstrat}  The starting field for the diophantine statement $D$ here is either $\bQ$ (or some fixed field $K$ containing $\bQ$) or a finite field $\bF_q$ (usually $\bF_p$).  To simplify we assume that we start either with $\bQ$ (or with coefficients in $\bZ$) or with $\bF_p$. Here is the procedure abstractly, starting with the definition of Galois stratification.  

In the rest of this subsection our goal is consider residue class fields. Then, \S\ref{quoteFrJ04} generalizes, as done in \cite{FrJ04}, but with different emphasis, how that gives   objects over $\bQ$. 

\begin{edesc} \label{gsconds} \item \label{gscondsa} There exists a stratification $\sS_{k}$ (a disjoint union of normal algebraic sets) $$\dot{{\cup}}_{j=1}^{l_{k}}X_{j,k},\text{ covering  }\afA^{m\np N_{k}}.$$
\item \label{gscondsb} For each $(j,k)$, $1\le j \le \ell_k$, there is a Galois (normal) cover $\hat \phi_{j,k}: \hat Y_{j,k}\to X_{j,k}$\footnote{By putting a $\hat {}$ over $\phi$ we are reminding that this is a Galois cover with automorphisms defined over the field of definition of the cover.}  with group $G_{j,k}$ and a collection, $\bfC_{j,k}$, of  conjugacy classes of $G_{j,k}$.\footnote{For the problems in this paper it suffices to take {\sl Frobenius classes}.  Each is a union of conjugacy classes of $G_{j,k}$ where if $g$ is in one of these, then so is $g^u$ if $(u,\ord(g))=1$.}
 
\item \label{gscondsc} Then, in place of $(x',\by_1',\by_2')\in (\text{or} \not\in) X(\bF_{q})$ we ask this. For $(\bx',\row {\by'} k)\in \sS_{k}(\bF_q)$,   
$$\text{if $(\bx',\row {\by'} k)\in X_{j,k}$, then the Frobenius,  $\Fr_{\bx',\row {\by'} k}\in G_{j,k}$ at the point, is in $\bfC_{j,k}$.}$$  
\end{edesc} 

\display{Ambiguity in the coefficients of equations} In \cite{FrJ04} the covers are always taken to be unramified, but the essential is whether having the Frobenius in the conjugacy classes is unambiguous. We have purposely left ambiguous  the coefficients of the equations defining the ingredients of \eql{gsconds}{gscondsa} and \eql{gsconds}{gscondsb}. Indeed, the procedure works like this. 

\begin{edesc} \label{goal} \item  \label{goala} G$_{\rm aap}$: The original statement is over $\bQ$ (or even over $\bZ$), but we don't care, for finitely many $p$, if it is not true or even for \eql{gsconds}{gscondsc} applicable. 
\item  \label{goalb} G$_{\rm allp}$: The original statement is over $\bZ$, meaningful for each prime $p$ and we have the option of deciding if, after finding it true for a.~a. $p$, whether it is true for all $p$. 
\end{edesc} 

\begin{defn} \label{gscharacteristics} For each Galois cover $\hat \phi: \hat Y\to X$ in a Galois stratification $\sS$, we need an affine space description of the underlying spaces, and of the group of the cover. Running over all possible $(j,k)$ we refer to this collection as the {\sl characteristics\/} of $\sS$. \end{defn} 

\hskip .25 in It is easiest to describe the characteristics when the spaces are described as explicit open subsets of affine hypersurfaces, as they are in \cite{FrJ04}. There are two possible situations where we might apply the stratification procedure assuming our initial goal is to decide if we can interpret it over $\bF_p$  by reduction mod $p$, and then if it is true for almost all $p$. 

These situations are compatible with Main Theorem \eqref{mainthm}. It uses coefficients initially over $\bQ$, producing a finite set, $M_{k\nm1}$,  of primes for which \eqref{excstratprimes} holds.   
\begin{edesc} \label{excstratprimes} \item \label{excstratprimesa} Either  the resulting (over) $\bQ$ equations don't make sense of \eql{gsconds}{gscondsc}; or the fibers of the underlying stratification pieces are not flat over $\Spec(\bZ_p)$; or 
\item  \label{excstratprimesb} the Chebotarev estimates of  \eqref{Frobest}  for some $q$ divisible by $p$ don't assure hitting all appropriate conjugacy classes. 
\end{edesc} 

The {\sl flatness statement\/} means that reduction mod $p$ preserves the characteristics of the stratification. In applying Denef-L\"oeser \S\ref{unifp}, what we  call the characteristics is strengthened, but the idea is the same. 

\begin{defn} \label{excstratprimesdef} In the $i$th step of the induction procedure, refer to the respective primes, the {\sl exceptional stratification primes} of \eqref{excstratprimes},  as $M_{i,a}$ and $M_{i,b}$, $i=k\nm1,\dots,0$, and their union by $M_i$. We include in $M_i$ the primes from $M_{i\np 1}$. So $M_0\supset M_1\supset \cdots \supset M_{k\nm1}$. \end{defn}  

It is for the strong Poincar\'e series result of Thm.~\ref{poinrat} that we have divided these considerations of primes. Following the stratification procedure format in Thm.~\ref{mainthm}, \S \ref{nonregcheb} gets into these significant issues: Reduction mod $p$; and the non-regular Chebotarev analog. 

The extra point, not discussed in \cite[Chap.~31]{FrJ04}, is how the diophantine illuminating properties of the Poincar\'e series  in Thm.~\ref{poincare} vary  as a function of $p\not \in M_i$. Especially, how they depend on the non-regular analog as illustrated by the \S\ref{schurcover} example.

\display{Inductive stratification notation} Below some additional notation replaces $k$ by the index $i$. For example, we replace  $N_k$ by $N_i=\sum_{j=1}^i n_j$. 
 Assume the problem has one quantifier for each block of variables. Our original example, illustrating a possible answer  in   \S\ref{Felgner} for Felgner, used two quantifiers. 
 
Galois stratification  allows eliminating one quantified block of variables, at a time, using a general,  enhanced, {\sl Non-regular Chebotarev Density Theorem}.  \S\ref{nonregcheb} explains how that contributes to  Main Theorem \eqref{mainthm} by which the elimination inductively  forms a sequence of Galois stratifications: $\sS_k, \dots \sS_0$, with triples $(X_{\bullet,u},G_{\bullet,u},\bfC_{\bullet,u})$, $u=k,\dots, 0$. 

For each $u$, the $\bullet$ indicates a sequence for $1\le j\le \ell_u$ of underlying spaces in the stratification that satisfy the following properties. 

As in \eql{gsconds}{gscondsb}, saying  $\Fr_{\bx',\row {\by'} {i}}\in \bfC_{\bullet,i}$ means that if the subscript point is in $ X_{j,i}$, then $\Fr_{\bx',\row {\by'} {i}}$ refers to  the conjugacy class value in the Galois cover  $\hat C_{j,i}\to X_{j,i}$. 

The variables $\bx$ of the Galois stratification $(X_{\bullet,0},G_{\bullet,0},\bfC_{\bullet,0})$ of $\afA^m$ are unquantified. For the  $i$th stratification $\sS_i$, the quantifiers are $\row Q i$;  we suppress their notation. 
\begin{defn} For $\bx'\in \bF_q^m$, we say $\sS_i(\bx';q)$  holds   
\begin{equation} \text{ if }\Fr_{\bx',\row {\by'} {i}}\in \bfC_{\bullet,i}\text{ for   }Q_1 \by'_1 \cdots Q_i \by'_i \in \afA(\bF_q)^{N_{i}}.\end{equation}  
A relative version concentrates at the $i$th position.  For $\bx',\by'_1,\dots, \by_{i\nm1}'\in \bF_q^{m\np N_{i\nm 1}}$    
\begin{equation} \text{ we say $\sS_i(\bx',\by'_1,\dots, \by_{i\nm1}';q)$ holds if } \Fr_{\bx',\row {\by'} {i}}\in \bfC_{\bullet,i}\text{ for   }Q_i \by'_i \in \afA(\bF_q)^{n_{i}}.\end{equation} \end{defn} 

The notation $\bar \sS_p\eqdef\{\sS_{p,k}, \dots, \sS_{p,0}\}$ from \eqref{stratproc}, for a sequence of Galois stratifications  -- with its attendent sequence of exceptional stratification primes $M_{k}\dots, M_0$ -- appears at the end of the Thm.~\ref{mainthm} statement.  Also, Thm.~\ref{mainthm} refers to the Frobenius in a finite Galois extension $\hat K/\bQ$, and the Frobenius (conjugacy class), $\Fr_p$  attached to a prime $p$ in  this extension. 

\begin{thm}[Main Theorem \cite{FrS76}]  \label{mainthm} 
Assume quantifiers $\row Q k$ and, $\sS_{k,p}$, an initial Galois stratification over $\bQ$, in situation \eql{goal}{goala}. Then, we may effectively form a stratification sequence $\bar\sS_p$, as above, for each $p\not \in M_0$ with the following properties. 
\begin{edesc} \item {\sl Local property}: Given $\bx',\by'_1,\dots, \by_{i\nm1}'\in \bF_p^{m\np N_{i\nm 1}}$, $\Fr_{\bx',\by'_1,\dots, \by_{i\nm1}} \in \bfC_{\bullet,i\nm1}$ if and only if $\sS_i(\bx',\by'_1,\dots, \by_{i\nm1}';p)$ holds. 
\item {\sl Inductive Result}: $\bx'\in \bF_p^m$, then $\Fr_{\bx'}\in \bfC_0$ if and only if $S_k(\bx',p)$ holds. \end{edesc}
Further, there is a finite Galois extension $\hat K/\bQ$ with group $G\eqdef G_{\sS_k}$ and a finite collection of stratification sequences ${}_j\bar \sS, j=1,\dots, \mu$, over $\hat K_p$ (see \S\ref{quoteFrJ04}) with this property. 
For each $p\not\in M_0$, $\bar \sS_p$ is the reduction of ${}_j\bar \sS$ mod $p$ for some $j$ depending only on $\Fr_p$ in $G$.  \end{thm} 

\begin{rem} The last paragraph in Thm.~\ref{mainthm} is there to compare the stratification procedure for a given $p$ to a process that works uniformly in $p$. We are giving several approaches: the Frobenius field approach (\S\ref{quoteFrJ04});  the finite field approach, and sometimes just deciding the truth of a diophantine statement for $p$. Still,  all of them -- including the ultimate count, as in the abstract, in the Poincar\'e series -- depend on locating absolutely irreducible varieties related to the stratification procedure.\end{rem} 

\begin{rem}[Extend Thm.~\ref{mainthm}] \label{extendq} For a given $p$, we may replace the quantification of the variables in $\bF_p$ by  quantification of its variables in $\bF_q$ with $q$ any power of $p$.  Again, the correct ${}_j\bar \sS$ depends only on $\Fr_q\in G$. See Rem.~\ref{largeq} for the adjustment that applies to considering  powers of $q$ divisible by $p\in M_0$. \end{rem} 

\begin{rem} \label{largeq} Notation of Def.~\ref{excstratprimesdef} gives a procedure starting with the subscript $i=k$, ending with $i=0$.  
Even, however, for primes in $M_0$ under hypothesis G$_{\rm allp}$ \eql{goal}{goalb}, the equation coefficients are in $\bF_p$. We may then perform the stratification procedure. That allows deciding for each such $p$ if the statement is true over $\bF_q$ for {\sl almost all\/} powers $q$ of $p$, so  long as $q$ is large enough for the Chebotarev conditions to hold. This extends Rem.~\ref{extendq}.  We cannot, however, expect the stratification to come from an object over $\bQ$ (by reduction mod $p$) or to be related in any obvious way to one of the stratifications ${}_j\bar \sS$ given in Thm.~\ref{mainthm}. \end{rem}

\subsection{The General Stratification procedure} \label{quoteFrJ04}  The qualitative point of the \NRC\ is this.  
\begin{equation} \label{NRCvalue} \begin{array}{c} \text{To know, for a Galois cover $\hat \phi:\hat Y\to X$ over a finite field, what conjugacy classes}\\ \text{ (or cyclic subgroups) generate decomposition groups running over $\bx\in X(\bF_q)$.} \end{array} \end{equation} 

\S\ref{nonregcheb}  shows the role of the Chebotarev Theorem  in Thm.~\ref{mainthm}. Especially that it gives a rough count of the $\bx' \in \bF_q^m$ for which $\sS_k(\bx',q)$ holds. This section outlines how \cite[Chaps.~30-31]{FrJ04} gives Thm.~\ref{mainthm}. We trace the general form of the stratification procedure. 

\begin{defn} A profinite group $G$ has the {\sl embedding property\/} if given covers\footnote{By a cover $\psi: H \to G$ of groups we mean an onto homomorphism.} of groups $\psi_{G,A}: G\to A$ and $\psi_{B,A}: B\to A$, for which  {\sl $B$ is a quotient of $G$}, then 
$$\text{there is a cover $\psi_{G,B}: G\to B$ with $\psi_{G,A}\circ \psi_{G,B}=\psi_{G,A}$.}$$  The name {\sl superprojective\/}  is thus apt for a group that is projective -- in the category of profinite groups -- with the embedding property. \end{defn} 

\hskip .25 in  That is, the embedding property is just projectivity with the extra condition that $\psi_{G,B}$ is a cover, given the necessary condition on $B$. This definition  starts the analogy between finite fields and the very large class of Frobenius fields -- again a reason for the existence of \cite{FrJ04} -- and whereby progress was achieved using Galois stratification \cite{FrHJ84}. 

\cite[\S24.1]{FrJ04} develops the theory of {\sl Frobenius fields\/} $M$:  
\begin{edesc} \label{frfi} \item  \label{frfia}  $M$ is a \PAC\ field (all absolutely irreducible varieties over $M$ have $M$ points), and;  \item  its absolute Galois group, $G_M$, has the embedding property. \end{edesc} $$\text{If $M$ is Frobenius, then $G_M$ is {\sl superprojective\/} \cite[Prop.~24.1.5]{FrJ04}.}$$ 

\display{The point of Frobenius fields} \cite[Prop.~24.1.4]{FrJ04}  is the replacement for Frobenius fields of the Chebotarev analog. This is applied to the Galois covers $\hat \phi: \hat Y\to X$ with group $G_{\hat \phi}$ that appear in a Galois stratification. 
This assumes $\hat Y$ is absolutely irreducible. Equivalently, $\hat \phi$ is a {\sl regular\/} cover, over $M$. 

The conclusion of this Chebotarev analog is that we know precisely what {\sl decomposition groups\/} --  among the $H\le G_{\hat \phi}$ -- that $\hat \phi$ has over fibers of $X(M)$. They are:  
$$\text{All subgroups $H$ that are quotients of  the absolute Galois group $G_M$. }$$  

The remainder of Chap. 30 through \S 30.5 is details of the Galois stratification procedure applied to the theory of {\sl all\/}  Frobenius fields containing a fixed given field \cite[Thm.~30.6.1]{FrJ04}. For simplicity assume we start over $\bQ$. Much of what looks technical is the pairing of decomposition groups of a cover $\hat \phi_i$  of $\sS_i$ with the decompositions groups that extend these to a given cover $\hat \phi_{i\np1}$ of $\sS_{i\np1}$ above $\hat \phi_i$.   

For example, \cite[conditions (2) of Lem.~30.2.1]{FrJ04} has this situation: $n_{i\np1}=1$ and $\by_{i\np1}=y_{i\np1}$. We continue Ex.~\ref{coverpair} in \S\ref{schurcover}.  Again, $M$ is notation for a Frobenius field. 

\begin{exmpl}[Pairing two covers] \label{coverpair} Suppose $X_2$ is a hypersurface in $\afA^{m\np1}$ defined by the equation $\{(\bx,y)\mid f(\bx,y)=0\}$ covered by $\hat \phi_2: \hat C_2\to X_2$. Also, that $\proj_\bx: \afA^{m\np1} \to  \afA^m$ restricted to $X_2$ is generically onto. So, it is a cover -- not likely Galois in a practical case -- over an open set in the image. Assume $f$ is absolutely irreducible of degree $u$ in $y$. 

Take a Zariski open set $X_1$ in $\afA^m$ and a Galois cover $\hat \phi_1: \hat C_1\to X_1$ having an open map $\phi: \hat C_1\to \hat C_2 \to X_2$ factoring through $X_1$. To satisfy other constraints in \cite{FrJ04}, like being unramified, it may not be surjective to $X_2$. Note: $G_{\hat \phi_1}$ has a natural degree $u$ permutation representation $T: G_{\hat \phi_1} \to S_u$ from its birational factorization through $\hat \phi_2$. 

The essence: Take for $\hat \phi_1$ the minimal Galois closure of the cover  $\proj_\bx\circ \hat \phi_2$.  Here is what we want from the elimination of quantifiers, if say $y$ was quantified by $\exists$. 

Given a Frobenius field $M$, consider subgroups $H_1\le G_{\hat \phi_1}$ that are decomposition groups for $\hat \phi_1$ over some $\bx'\in X_1(M)$ that also fix a point $(\bx',y')\in X_2(M)$. That is, $H_1$ fixes a letter in the representation $T$, so it restricts to a decomposition group $H_2'$ of $\hat \phi_2$. Then, consider $\sH$, the conjugacy classes of groups $H_2\le G_{\hat \phi_2}$ in Galois stratification data already attached to this cover in the Frobenius field case. 

How to decide --  for  eliminating the quantifier $\exists$ --  if $H_1$ will be in the conjugacy classes attached to the cover $\hat \phi_1$? It is, if the $H_2'$ described above corresponding to $H_1$, is in $\sH$. \end{exmpl} 

\begin{rem}[Finding $H_1$ in Ex.~\ref{coverpair}] \label{fieldcrossing} The last line on finding those $H_1\,$s is what the Chebotarev analog quoted above does. In detail, you construct an absolutely irreducible variety for which an $M$ point corresponds to such an $H_1$. Here I allude to Chebotarev's {\sl field crossing argument\/} -- which is all over \cite{FrJ04} starting on p.~107 -- and again in \S\ref{choices} and \S\ref{schurcover}. \end{rem}  

\display{Continuing the stratification procedure}\  Reducing the general situation to Ex.~\ref{coverpair} consists of many, conceptually easy, steps. A large collection of Frobenius fields satisfies the general theory of Frobenius fields containing any Hilbertian field $K$. 

\hskip.25 in To compare with finite fields \cite[\S 30.6]{FrJ04} cuts things back from all Frobenius fields to much smaller collections of Frobenius fields. It does so by starting with $\sC$ some {\sl full formation\/} of finite groups. Then, it limits $G_M$ to be in the pro-$\sC$ groups,   

Cutting down further, \cite[p.~720]{FrJ04} considers a fixed superprojective group $U$. It restricts to those $M$ containing $K$ with $G_M$ in the  collection of quotients of $U$. Finally, it goes to the case totally compatible with finite fields: $U$ is the profinite completion of $\bZ$. 

Restricting the full theory of Frobenius fields containing $K$ to these special cases is quite conceptual. The final realization is that in this last case the stratification procedure is done entirely over $\bQ$. Now reduce the whole stratification apparatus mod $p$ for most primes $p$. This amounts to dealing with covers, over $\bQ$ (or later, $\bF_q$) especially those that aren't regular. 

Thm.~\ref{nonregcover}  \cite[Lem.~1]{Fr74} handles this.  For a Galois cover $\hat \phi: \hat Y\to X$ (again take both $X$ and $\hat Y$ to be normal) with $X$ absolutely irreducible having group $\hat G$ \cite[Lem1]{Fr74}.\footnote{This lemma does do the curve case, but the details necessary to generalize are the elimination results typical in commutative algebra.} \S\ref{schurcover}  shows why this non-regularity, not easily eliminated, occurs. Take $\hat \phi$ to be a $K$ irreducible cover, $K$ a number field, with ring of integers $\sO_K$. 

Assume $\hat K$ is the minimal field over which $\hat \phi$ breaks into absolutely irreducible (Galois) covers of $X$. Take any component $\phi': Y'\to X$. It will be Galois with group identified with $G'\eqdef \{g\in \hat G\mid g \text{ maps } Y'\to Y'\}$. The Galois extension $\hat K/K$ has group $\hat G/G'$. For $\bp$ a prime of $\sO_K$ unramified in $\hat K$, denote an element of the  conjugacy class of the Frobenius of $\bp$ by $\Fr_\bp$, and its fixed field (resp.~residue class field)  in $\hat K$ by $\hat K_{\bp}$ (resp.~$\bF_\bp$). Then, let $\hat \phi_\bp: \hat Y_\bp\to X$ be the union of the conjugates of $\phi': Y'\to X$ by $\Fr_\bp$. 

\begin{thm} \label{nonregcover}   Excluding a finite set $B$ of computable primes $\bp$ of $\sO_K$, reduction of $\hat \phi$ mod $\bp$ has an $\bF_{\bp}$ Galois cover component isomorphic to $\hat \phi_\bp$ by the reduction map. Take the conjugacy classes, $\bfC_\bp$ associated to $\hat \phi_\bp$ to be those from $\bfC$ whose restriction to $\hat K_{\bp}$ act trivially. 
\end{thm} 

\begin{rem} The classes $\bfC_\bp$ that appear in Thm.~\ref{nonregcover} are precisely those that can be realized as Frobenius elements,  as in Comment on \eql{chebuses}{chebusesb} in \S\ref{othcheb}.  \end{rem} 

\begin{rem}[Proofsheet changes, Edition 2 vs Edition 3 of \cite{FrJ04}] \label{corrections} We list the two that concern the topic of this paper. In both editions a key stratification lemma repeats on consecutive pages as Lem.~30.2.6 and Lem.~30.3.1 (including the material in the 2nd edition that starts at the bottom of p.~711). Also, the references on p. 760 in Edition 2 between [W.D. Geyer and M.~Jarden] and [W.~Kimmerle, R. Lyons, R. Sanding and D.N. Teaque] are missing in Edition 2, but appear in Edition 3 on pgs.~766--771. \end{rem}  

\subsection{Choices and  the Non-regular Chebotarev}  \label{nonregcheb}  From the view of deciding diophantine statements, \cite{FrS76} is driven by actual diophantine applications, while \cite{FrJ04} is stronger on logic's  theories of fields. This paper returns to the first viewpoint: Galois stratification as  enhancing understanding specific problems. 

The variance of one diophantine problem with the prime $p$ starts with using the Weil estimate systematically in \S\ref{choices}, including finishing off Felgner in \S \ref{concfelgner}. Variance gets enhancement from attaching to the problems a Poincar\'e series, the topic of \S\ref{diophinvar}. 

\subsubsection{Choices in stratifying} \label{choices} 
Being  courser -- not stratifying excessively -- on the formation of Galois stratifications allows staying closer to classical problemsand explicit computation.  At that, the original paper \cite{FrS76} and \cite{FrJ04} are at opposite ends of the spectrum, despite the latter's details. 

The former suggests restricting to flat covers when appropriate, involving blocks of variables of maximal length when possible. The latter takes blocks with just one variable and insists on unramified covers given by $\Spec(S)\to \Spec(R)$ with $S$ generated, as a polynomial ring, over $R$, by a single element having no discriminant.\footnote{This guarantees that a Frobenius element is a well-defined conjugacy class, without demanding extra conditions on $\bfC$. In, however,  classical problem settings, using flexibility on $\bfC$ may be a good idea.}  \S\ref{quoteFrJ04} outlined using the finer stratification procedure, as does this section which relies only on the original Weil estimate. 

On the other hand, \S\ref{schurcover} gives our main example, a general situation that concludes with a Galois stratification containing just one cover, and by which we painlessly see  infinitely many of the Poincar\'e series coefficients. 

Now we show how the Chebotarev analog works, with a corollary to Main Theorem \ref{mainthm} that counts $\bx'\in \bF_q^m$ for which $S_k(\bx',q)$ holds (for $q\not \in M_0$). First, consider going from $\sS_{i\np1}$ to $\sS_{i}$, dealing with the restriction of the stratification for $\sS_{i\np1}$ along the fibers of the \begin{equation} \label{progitoi-1} \text{natural projection  $ \afA^{m\np N_{i+1}}\ \longmapright{\pr_{i\np 1,i}} {40}\  \afA^{m\np N_{i}}$. }\end{equation}   

With $\ell_0$ the number of elements in $\sS_0$, and each $1\le j\le \ell_0$ this gives the following  \cite[Thm.~p. 104]{FrHJ94}:
\footnote{Qualitatively \cite[Prop.~2]{Fr74} sufficed for \cite{FrS76}, but there are more details in the unspecified constants, especially the dependency on the characteristics of the covers, in this reference.} 
\begin{edesc} \label{Frobest} \item an integer $r_j$ between 0 and $N_{k}$ and $\mu_j\in \bQ^+$, a function of elements in $\bfC_{j,0}$;  
\item \label{Frobestb}    giving the count, $B_j$,  of $\bx'\in \afA^m$   with $\bx'\in X_{j,0}$ and $\Fr_{\bx'}\in \bfC_{j,0}$. 
$$\text{Either: \ $(\dagger)\ B_j$ is 0;   {\sl or\/} $(\dagger\dagger)\   |B_j -\mu_j q^{r_j}|/ q^{r_j\nm 1/2}$ is bounded in $q$.}$$ 
\end{edesc} 

These estimates  \cite[\S3, especially Lem.~3.1]{FrHJ94} are based on \cite{LW54}. They are explicit, much better than {\sl primitive recursive\/}, as in  the application  concluding \cite{FrS76}.   In \cite[Thm.~6.1, 6.3 et. al.]{FrHJ94}, the main ingredients is this. Apply a Chebotarev Density Theorem version (as in \cite[Prop.~2]{Fr74}) to a pair consisting of the Galois cover $\hat \phi: \hat C\to X$ to estimate the number of points $\bx'\in X(\bF_q)$  for which $\Fr_{\bx'}$ is in the conjugacy classes $\bfC$ attached to $\hat \phi$. Comments on \eql{chebuses}{chebusesa} show why, in the inductive procedure, we are producing new covers with new  Galois closures and corresponding conjugacy classes. 

Tesselating $X$ with hyperplane sections -- akin to \cite{LW54}  -- and using Chebotarev's own field crossing argument, reverts this to Weil's Theorem on projective nonsingular curves over finite fields. \cite[\S4]{FrHJ94} explicitly  traces these classical results. Yet, it still has an error estimate. So, it doesn't imply the rational function Poincar\'e series results  in \S\ref{poincare} or \S\ref{unifp}. 

\subsubsection{Other Chebotarev Points} \label{othcheb} 
\begin{edesc} \label{chebuses}  \item \label{chebusesa} The actual quantifier elimination moving from $\sS_{i\np 1}$ and $\sS_{i}$ also uses \eqref{Frobest}. Indeed, that shows why there is no elimination of quantifiers through {\sl elementary\/} statements. 
\item \label{chebusesb} Possibilities $(\dagger)$ and $(\dagger\dagger)$ in \eql{Frobest}{Frobestb} give very different error estimate contributions. 
\item \label{chebusesc} The distinctions in \eql{chebuses}{chebusesb} arise in a host of problems as illustrated in \S\ref{schurcover}. \end{edesc} 
 
\display{Comment on \eql{chebuses}{chebusesa}}  Consider restriction of one of the terms of the stratification of $\sS_{i\np1}$ to the fibers $\pr_{i\np1,i}$ of the projection \eqref{progitoi-1}. Much of the \cite{FrS76} proof assures the elimination theory allows applying \eqref{Frobest}, so as to pick out the conjugacy classes that will appear in each of the terms of the $\sS_{i}$ stratification. \S\ref{quoteFrJ04} has that, though using Frobenius fields there avoided reference, at that point, to the Weil result. 

So, for the quantifier $\exists$, the Chebotarev analog is essentially to assure that any conjugacy class that should occur in $\sS_{i}$ actually  does. Likewise, for $\forall$, that no conjugacy class that could change the result of the problem is excluded. Nevertheless, Chebotarev is giving error estimates, and not {\sl precise\/} values that contribute to refined invariants as  in \S\ref{diophinvar}. 

\display{Comment on \eql{chebuses}{chebusesb} }  The distinction between $(\dagger)$ and $(\dagger\dagger)$ is consequent on the adjective {\sl non-regular\/} analog of the Chebotarev (any dimensional base) version. Before  \cite[\S2]{Fr74} it had been traditional to make an unwarranted assumption. That is, when considering a cover $\phi: C\to X$ -- say of normal varieties -- defined and absolutely irreducible, say, over a field $K$, that some kind of manipulation would allow assuming that the functorially defined Galois closure $\hat \phi: \hat C\to X$ of the cover could also be taken over $K$. 

That won't work in considering the possibilities, based on one Galois stratification, in varying $q$, even in residue classes of a number field. Here is why (eschewing cautious details). Suppose a component of the Galois closure cover, $\hat\phi$, has definition field $\hat K\not = K$ (with $K$ a perfect field). Assume $K$ is a number field with ring of integers $\sO_K$. 

Then, as we vary the residue class field $R_\bp=\sO_K/\bp$, the corresponding Frobenius $\Fr_{\bx'}$ for $\bx'\in \bF_q^m$ must restrict to the Frobenius $\Fr_{\bp}$ on the residue class field $R_\bp$. Suppose, for example, the order of $\Fr_{\bp}$  does not divide the order of an element $g\in \bfC_{j,0}$.  

Then,  that conjugacy class of $g$ cannot  possibly be  $\Fr_{\bx'}$.\footnote{That is, there is no error estimate for non-achievement of that class as a Frobenius, say by a bounding constant; it is just not achieved.}  That, however, is the only obstruction by the general Chebotarev result -- the meaning in \eql{excstratprimes}{excstratprimesb} of hitting the correct classes -- for realizing an element of $\bfC_{j,0}$ as a Frobenius, so long as $q$ is sufficiently large. 

\display{Comment on \eql{chebuses}{chebusesc}} The comment on \eql{chebuses}{chebusesb} alludes, for $R=\sO_K$, to the $(\dagger)$ and $(\dagger\dagger)$  conditions varying with the residue class field $R_\bp=\sO_K/\bp$ in actual problems. Further, the value of the Frobenius in an extension of $\hat K/K$ attached to the characteristic of $R_\bp$ measured this. The analog is true for problems over a given finite field: The Frobenius in an extension $\hat \bF_q$ attached to the problem measures the variance with changing the extension of $\bF_q$. Both justify the significance, and resistence to elimination, of the word {\sl non-regular\/} in the Chebotarev analog. \S\ref{schurcover} is an example that makes this explicit. 

\subsubsection{Main Theorem \ref{mainthm} implies \lq\lq No!\rq\rq\  to Felgner}  \label{concfelgner} 

Back to Felgner's Question with  $m=1$ and running over elements of $\bF_{q^2}$ (not $\bF_q$). Then, the elimination of quantifiers has reduced Felgner's question to this.

\begin{edesc} \label{counting} \item \label{countinga}  Show the following is impossible for $q$ large:  $$\text{with $M_w=|x'\in X_{w,0}(\bF_{q^2})\text{ and } \Fr_{x'}\in \bfC_{w,0}|$,  $M_{\sS_0}\eqdef \sum_{w=1}^{\ell_0}M_w =|\bF_q|=q$.}$$ 
\item Main Theorem \ref{mainthm} and \eqref{Frobest} implies  for each $w$,  $M_w$ is either bounded (in $q$)  or asymptotic to  $t_wq^2$ for some nonzero $t_w$. 
\item Neither is a bounded distance from $q$. Therefore,  no elementary formula distinguishes $\bF_q$ within  $\bF_{q^2}$, so long as $q$ is large. \cite[\S0]{FrHJ94}
\end{edesc}  

Indeed, the same argument shows there is no need to take $m=1$. No matter what formula, no matter the number of variables included in $\bx$, you can't get  $q$ as the number of $\bx'$ counted by  \eql{counting}{countinga}, so long as $q$ is large.\footnote{\cite{FrHJ94} exists because the authors of \cite{CLMa92} insisted in their first version that Galois stratification couldn't handle Felgner's question. This  despite others at their conference, including one of the authors of \cite{DL98} -- I was not -- telling them that was wrong, much akin to the end of  \cite{FrS76}.}
So we used quantitative counting to {\sl exclude\/} the elements of $\bF_q$ as a result. This was rather than finding a qualitative device that eliminated existence of a formula that precisely nailed $\bF_q$ among the elements of $\bF_{q^2}$.\footnote{That the count is bounded by a constant  can't be excluded; see the {\sl exceptional covers\/}  of \S\ref{unifp}.}  

Is that  satisfactory? It gets to the heart of the  Poincar\'e series/zeta function approach, which primarily aims to count points satisfying equations.  
 
\section{Diophantine invariants} \label{diophinvar}  \cite[\S7.3]{Fr05} discusses the history of attaching a Poincar\'e series (and zeta function) to diophantine problems attached to Galois stratifications. We refer to  some of its highlights. This section pushes Chebotarev into the background; replacing it with expressions in precise point counts. \S\ref{poincare} gives the Poincar\'e series definition. 

\S\ref{questpoin} goes through the series' properties, especially Rationality Thm.~\ref{poinrat}.  These estimates use {\sl Dwork cohomology\/} and specific results of Bombieri applied to it. 

\S\ref{schurcover} ties together all threads of this paper with one specific problem -- having a vast practical literature -- that uses the Poincar\'e series. Then, \S\ref{chowmotives} briefly discusses the artistic extension -- of Denef and L\'oeser -- from Galois stratification to {\sl Chow Motive\/} coefficients. Thus, \S\ref{unifp} is in the service of enhancing uniformity in $p$.

\subsection{Poincar\'e series vs coefficient estimates} \label{poincare}  We have already seen that the Galois stratification procedure is not canonical. So, it makes sense to address whether there is a {\sl homotopy theoretic\/}  approach based on the category of Galois stratifications that clarifies a natural equivalence on stratification. That is what these Poincar\'e series test.\footnote{Take {\sl homotopy theoretic\/} to mean that an outcome expressable in terms of some kind of cohomology. One that equivalences among structures related to Galois stratification or variants that result in the same cohomology results. In a sense that is the point of motivic cohomology.} 

\eqref{poinform} is the definition of the Poincar\'e series attached to the Galois stratification $\sS_k$ over $\bF_q$ with triples $(X_{\bullet,k},G_{\bullet,k},\bfC_{\bullet,k})$, and its quantifiers $\row Q k$. Main Theorem \ref{mainthm} replaces this by quantifying those points of the last stratification term whose Frobenius values are inside the requisite conjugacy classes. Though more complicated than counting points on a variety over a finite field, it is sufficiently akin to naturally extend classical methods.  

The Poincar\'e series, in a variable $t$,  for a given $q$ attached to Galois Stratification $\sS_k$ has this form with the coefficients  $\mu(\sS_k,q, m)$ explained below: 
\begin{equation} \label{poinform}  P(\sS_k)_q(t) =  \sum_{m=1}^\infty  \mu(\sS_k,q, m) t^m.\footnote{The notation of \S\ref{hrushwan} shows why we didn't use the simpler notation $\mu(\sS_k,q^m)$.}\end{equation} 
Those $\mu(\sS_k,q, m)\,$s do depend on the quantifiers ${\pmb Q}$. If the Galois stratification was an ordinary elementary statement, then to simplify we would abuse the notation by placing them outside reference to the variables. We mean   each quantifier $Q_i$, respectively,  applies to the variables $\by_i$, $i=1,\dots k$. For example: ${\pmb Q}\, \bx\,  \by_1 \,\by_2$ is $\bx\, Q_1\by_1 Q_2\by_2$.   

So, similar to this, when quantifying the placement of the Frobenius elements, denote the quantified version by ${\pmb Q}  \Fr_{\bx',\row {\by'} k}\in \bfC_{j,k} $  
\begin{equation} \label{pointscounted} \text{ running over }(\bx',\row {\by'} k)\in X_{j,k}(\bF_{q^m}) (\bx',\row {\by'} k), 1\le j\le\ell_k. \end{equation} 
Then, the coefficient $\mu(\sS_k,q, m)$ in \eqref{poinform} is the point count of  those $\bx'\in X_{j,k}(\bF_{q^m})$ for which ${\pmb Q}  \Fr_{\bx',\row {\by'} k}\in \bfC_{j,k} $  is constrainted to  running over the  quantified $y$-variables, $\row {\by'} k$, with values in $\bF_q$.   
Main Theorem \ref{mainthm}  allows replacing  $\mu(\sS_k,q,m)$ by $\mu(\sS_0,q,m)$:  counting (rather than estimating as in \eql{Frobest}{Frobestb}) those \begin{equation} \label{S0coef} \text{$\bx'\in X_{j,0}(\bF_{q^m})$ for which $\Fr_{\bx'}\in \bfC_{j,0}, 1\le j\le \ell_0$.} 
 \end{equation}    
 
 \subsection{Poincar\'e properties} \label{questpoin}  
 
 \cite[Chap.~30 and 31]{FrJ04} (Chaps.~25 and 26 in the 1986 edition, pretty much the same)  have  complete details on the Poincar\'e and Zeta properties. 
A {\sl Zeta function\/}, $Z(t)$, has an attached {\sl Poincar\'e series\/} $\tP(t)$. This is given by the logarithmic derivative: $$t\frac{ d}{dt}\log(Z(t))=\tP(t).$$ Add that $Z(0)=1$, and each determines the other.  The catch: $Z(t)$ rational (as a function of $t$) implies $\tP(t)$ rational, but not always the converse. 

Suppose $D$ is an elementary diophantine statement, with quantifiers as we started this paper. Then, as earlier,  take $D$ in place of $\sS_k$ in $P(\sS_k)_q(t) =  \sum_{m=1}^\infty  \mu(\sS_k,q, m) t^m$, and consider the coefficients referencing $(D,{\pmb Q}, m)$ in place of  $\mu(\sS_k,q, m)$. 

I don't know when Ax introduced considering such  coefficients. He suggested to me meaningfully computing them at IAS in Spring '68. Originally,  I introduced the Galois stratification procedure to do just that, and to conclude the following result. Again, suppress notation ${\pmb Q}$ for the quantifiers and give the result explicitly just when $R=\bZ$ in \eqref{tradff}. 

The adjustments are clear for when $R$ is a given finite field, or ring of integers of a number field.  Use the notation around Main Thm.~\ref{mainthm}. 

\begin{thm}[Poincar\'e rationality] \label{poinrat} 
For each prime $p\not \in M_0$, $P(\sS_k)_p(t)$   is a rational function $\frac{n_p(t)} {d_p(t)}$, with $n_p,d_p\in \bQ[t]$ and computable.  The corresponding $Z(\sS_k)_q(t)$ has the form $\exp(m_p^*(t))(\frac{n_p^*(t)} {d_p^*(t)})^{\frac 1 u_p}$ with $m_p^*,n_p^*,d_p^*\in \bQ[t]$ and $u_p\in \bZ^+$ computable. Further, there are bounds, independent of $p$, for the degrees of all those functions of $t$. For any particular prime $p$ all functions are computable (see \S\ref{chowmotives}). 
  \end{thm}

\display{Comments on the proof of Thm.~\ref{poinrat}}  We start with highlights from \cite[\S31.3]{FrJ04}  (or 1986 edition, \S26.3;   essentially identical) titled: Near rationality of the Zeta function of a Galois formula. A similar result bounds the degrees even if $p\in M_0$, assuming \eql{goal}{goalb}. 

\hskip .25in   As, however, Rem.~\ref{largeq} notes, there won't be any expected relation with results for $p\not\in M_0$. We cannot use the uniform estimate on the characteristics of the Galois stratification since they don't come from a reduction of a uniform stratification object as given in Thm.~\ref{nonregcover}.  Refer to the stratifications for $p\in M_0$ as {\sl incidental}, with their incidential estimates. What we say here applies equally to uniform and incidental stratifications. 

The conclusion of the Galois stratification procedure over the $\bx$-space gives this computation for \eqref{S0coef}. Sum the number of $\bx$ with values in $\bF_{p^k}$ for which the Frobenius falls in conjugacy classes attached to the stratification piece  going through $\bx$.  

The expression of that sum in {\sl Dwork cohomology\/} is what makes the effectiveness statement in the Thm.~\ref{poinrat} possible. That also suggests its direct relation to Denef-Loeser. An ingredient for that is a formula of E.~Artin.  It computes any function on a group $G$ that is constant on conjugacy classes as a $\bQ$ linear combination of characters induced from the identity on cyclic subgroups of $G$. 

\display{Additional historical comments} 
A function on $G$ that is 1 on a union of conjugacy classes, 0 off those conjugacy classes, is an example. \cite[p.~738-739]{FrJ04} (1986 edition p.~432-433) recognizes L-series attached to that function as a sum of L-series attached to those special induced characters. I learned this from  \cite[p.~222]{CaFr67} and had already used it in \cite[\S2]{Fr74}. Kiefe -- working with Ax -- learned it, as she used it in \cite{Ki76},  from me as a student  during my graduate course in Algebraic Number Theory at Stony Brook in 1971. The core of the course were notes  from Brumer's Fall 1965 course at UM.\footnote{I submitted my paper in 1971. It had five different referees.} 

Kiefe, however, applied it to a list-all-G\"odel-numbered-proof procedure; not to Galois stratifications I  showed her (see my Math Review of her paper, Nov.~1977, p.~1454). Consider the  identity representation induced from a cyclic subgroup $\lrang{g}, g\in G$. This  L-Series is  the zeta function for the quotient cover by $\lrang{g}$ (exp. 7-9, p.~433, 1986 edition of \cite{FrJ04}).  

\display{We do use a zeta function} Given a rational function in $t$, its  {\sl total degree\/} is the sum of the numerator and denominator degrees; assuming those two are relatively prime. The 1986 edition, Lem.~26.13 refers to combining \cite{Dw66} and \cite{Bo78} to do the zeta function of the affine hypersurface case for  explicit bounds -- dependent only on the degree of the hypersurface -- on the total degree  of the rational functions that give these zeta functions. 

\hskip .1inThen, some devissage gets back to our case, given explicit computations dependent only on the degrees of the functions defining these algebraic sets.  Lem.~26.14 assures the stated polynomials in $t$ have  coefficients in $\bQ$, and it explicitly bounds their degrees. The trick  is to take the logarithmic derivative of the rational function. Then, the Poincar\'e series coefficients are power sums of the zeta-numerator zeros minus those  of the zeta-denominator zeros. Using allowable normalizations, once you've gone up to the coefficients of the total degree, you have determined the appropriate numerator and denominator of $\tP(t)$. 

One observation: We are left to  uniformly bound in $p$ the degrees of the zeta polynomials, etc. This follows from the comments above Thm.~\ref{nonregcover} giving the uniformity in $p$ in the characteristics of the reduction of the stratification. Apply this to  the degrees of polynomials describing the affine covers for Dwork-Bombieri.    

\subsection{Chow Motive Coefficients} \label{unifp}  As stated in \S\ref{choices}, the error estimate that allowed the elimination of quantifiers isn't appropriate for concluding either the estimate of the degrees of the polynomials in the Poincar\'e series, or that it is a rational function.  

\hskip .1in It is sensible to use as coefficients of the Poincar\'e series actual Galois stratifications. For $p\not\in M_0$ those coefficients sum over the Galois covers $\hat \phi: \hat Y \to X$ in the stratification $\sS_0$ those $\bx'\in X(\bF_q)$  for which $\Fr_{\bx'}\in \bfC$. It also works to replace the count on the stratification pieces with absolutely irreducible algebraic varieties that give the same counts. This uses the field crossing argument mentioned several times previously, as in Rem.~\ref{fieldcrossing}.  

The topic here is  enhancing  the Poincar\'e series coefficients, extending them to  {\sl Chow motive coefficients}; the last of the coefficient choices given in the abstract.  After the example of  \S\ref{schurcover}, \S\ref{chowmotives}  briefly discusses the abstract setup of \cite{DL02} and \cite{Ni07}. These have their own expositions on Galois Stratification. Further discussion of these will occur in the extended version of this paper. 

\subsubsection{Distinguishing special primes} \label{schurcover}  

We produce an example where the associated Poincar\'e series over $\bF_q$ has some coefficients that are {\sl polynomials\/} in $q$. Yet, other coefficients are functions of the Frobenius $\Fr_q$ evaluated on a Chow motive: a linear expression including the (nontrivial) $\ell$-adic cohomology of a nonsingular variety.  \cite[\S3]{Fr05} takes a general diophantine property and casts it as an umbrella over two seemingly distinct diophantine properties about general covers. 

These properties fit a {\sl birational\/} rubric, or what we call {\sl monodromy precision}. That is, the Galois closure of the covers alone, together with  the corresponding permutation representations attached to their geometric and arithmetic monodromy, guaranties, precisely, their defining diophantine property. The \NRC\ -- as always -- gives the count of achieved conjugacy classes {\sl roughly\/}. Here, though, there will be no error term -- unlike the Comment on \eql{chebuses}{chebusesb} in  \S\ref{othcheb} -- even though there is achievement of nontrivial conjugacy classes. 

That allows stating, should we start with such a cover over a number field $K$, what are the primes for which the diophantine property has a precise formulation over a given residue class field $\sO_K/\bp=\bF_\bp$, as in \cite[Def.~3.5 and Cor.~3.6]{Fr05}. Our example uses the simplest case: the {\sl exceptional cover\/} property. 

\begin{exmpl}[Ex.~\ref{coverpair} continued]  \label{coverpaircont} 
For this, in Ex.~\ref{coverpair}, take the following special case over $K$. The hypersurface  in $\afA^{m\np1}$ is still defined by an equation $X_2= \{(\bx,y)\mid f(\bx,y)=0\}$, with $f$ absolutely irreducible (over $K$) of degree $u>1$ in $y$. Now, though, the cover $\hat \phi_2$ is trivial (of degree 1). As before consider $\proj_\bx: \afA^{m\np1} \to  \afA^m$ restricted to $X_2$ and, these diophantine statements:   
\begin{equation} \label{expprop} \begin{array}{rcl} D_\bp(\bx):& \exists\ y \in \bF_\bp &\mid f(\bx,y)=0; \\ D_\bp:& \forall\  \bx\in \bF_\bp^m\  \exists\  y \in \bF_\bp&\mid f(\bx,y)=0; \\ 
D: & D_\bp \text{ is true for $\infty$-ly many $\bp$.} 
 \end{array} \end{equation}   

Continue with the Ex.~\ref{coverpair} notation, and the formation of the Galois cover $\hat \phi_1: \hat C_1\to X_1$ with $X_1$ Zariski open in $\afA^m$, and with group $G_{\hat \phi_1}$ having  its natural and faithful, transitive, degree $u$  permutation representation $T$.  Consider the projective normalization, $\hat C_1^\dagger$,  (resp.~$X_2^\dagger$) of $\hat \phi_1$ (resp.~$X_2$) in the function field of $\hat C_1$  (resp.~of $X_2$) $$\hat\phi_1^\dagger:  \hat C_1^\dagger \to X_2^\dagger \to \prP^m, \prP^m\supset \afA^m.\footnote{$\hat C_1^\dagger$ is the projective normalization of $\prP^m$ in the function field of  $\hat C_1$; a canonical process.}  $$ 
Now extend the diophantine expressions of \eqref{expprop} to include $X_2^\dagger$. For example, $D_\bp(\bx)$, for $\bx\in \prP^m(\bF_q)$ means $\exists$ a point of   $X_1^\dagger(\bF_q)$ above $\bx$. A similar meaning is given to $D_\bp$.\end{exmpl} 

We make a simplifying assumption in Ex.~\ref{coverpaircont} to match up with \S\ref{chowmotives}. \begin{equation} \label{nonsing} \text{Not only is $X_2$ normal, but $X_2^\dagger$ is nonsingular.}\end{equation}  Though the spaces are given by projective, not affine, coordinates, we form a single cover that we may regard as a Galois stratification $\sS_0$, with one stratification piece: $\hat \phi_1^\dagger: \hat C_1^\dagger \to \prP^m$ with attached conjugacy classes 
$$\bfC_1=\{g\in G_{\hat \phi_1}\mid T(g) \text{ fixes a letter in the representation}\}.$$ 

\begin{prop} \label{exceptprop} Assume \eqref{nonsing}. There is a finite set, $M_0$, of primes $\bp$ of $\sO_K$  for which, given $\bp\not \in M_0$, then $D_ \bp$ is true if and only the following equivalent conditions hold: 
\begin{edesc} \label{exccond} \item \label{excconda}  $\Fr_{\bx'}\in \bfC_1$ for each $\bx'\in \prP^m(\bF_\bp)$. 
\item \label{exccondb}  There are infinitely many $\bp$ for which the equivalence of \eql{exccond}{excconda} holds. 
\end{edesc} Assume \eqref{exccond}. Then for some finite set $M_0'\supset M_0$, for each $\bp\not \in M_0'$ for which $\eql{exccond}{excconda}$ holds,   
\begin{equation} \label{exccondbp} \begin{array}{c} \text{if $\bF_\bp=\bF_{q_0}$, \eql{exccond}{excconda} holds with $\bF_{q_0^m}$ replacing $\bF_\bp$. Then, for $\infty$-ly many $m$, } \\ \text{the $m$th coefficient of the Poincar\'e series $P(\sS_0)_\bp(t)$ for   $\sS_0$  is $\frac{q_0^{m\np1}-1}{q_0\nm1}$.}\end{array} \end{equation} We also note: 
\begin{equation} \label{exccondbpnot} \text{There are also infinitely many $\bp$ for which \eql{exccond}{excconda} does not hold. (See below.)} \end{equation} 
 \end{prop} 

\begin{defn}[Exceptional primes] The primes $\bp$ for which \eqref{exccondbp} holds are called the {\sl exceptional primes\/}, $E_{\hat \phi_1^\dagger}$, of $\hat \phi_1^\dagger$ (or of any other object natural attached to it). The Main point is there are infinitely many of them, and their Poincar\'e series gives them away. Such a cover is called an {\sl exceptional cover\/} over $K$.  \end{defn}

Note, however, this is a case of a regular cover whose Galois closure is {\sl not\/} regular. What this says in \eql{exccond}{exccondb} is that only the elements of $\bfC_1$ are achieved as Frobenius elements for the primes satisfying the conditions \eqref{exccond}. For example, for those primes $\bp$ (with $|\bp|$ suitably large) for which reduction $\mod \bp$ gives a regular function field extension, the Chebotarev Theorem says \eql{exccond}{excconda} will {\sl not\/} hold. It is very elementary -- an easy case of Chebotarev for number fields -- that there are infinitely many such $\bp$. 

\cite{DaLe63} considered this hyperelliptic pencil with parameter $\lambda$: $y^2 - f(x) + \lambda$, $f \in \bF_p(x)$. The difference between the expected number, $p$, of $\{(x,y)\in (\bZ/p)^2\}$ and the actual value, $V_{p,\lambda}$, given by Weil's theorem caused them draw the following conclusion. 
$$\begin{array}{l}\text{There is a constant $c_f>0$ such that: } \\ \text{Running over $\lambda\in \bF_p$, $\sum_{\lambda\in \bF_p}(p-V_{p,\lambda})^2> c_fp^2$, if and only if $f$ is not exceptional.}\end{array}$$ 

\cite{Ka88} (recounted in \cite[\S7]{Fr05}) used a generalization of the Davenport-Lewis error term to count the summed squares of the multiplicity of  (completely reducible) components of the monodromy (fundamental group) action on the 1st complex cohomology of a fiber of any family of nonsingular curves over a base $S$. His technique -- also coming upon the exceptionality condition -- used reduction mod $p$ to get to the results of \cite{De74}.  

\begin{rem} \label{langlands}  It has long been known that there are many exceptional covers. Between \cite[\S2]{Fr78}, \cite{GMS03} and \cite[\S6.1 and 6.2]{Fr05} those with $m=1$ and $f(x,y)=f(y)-x$, $f$ a rational function over a number field. The 1st and 3rd of these connect to Serre's Open Image Theorem. That paper also introduced natural zeta function test cases of the Langland's program. One problem is a standout.  For the wide class of these related to the $\GL_2$ case (and only these) of Serre's Open Image Theorem, the precise set  $E_{\hat \phi_1^\dagger}$ is a mystery appropriate for the nonabelian class field theory of the Langland's program \cite[\S6.3]{Fr05}. \end{rem}

\subsubsection{Denef-L\"oser and Chow motives}  \label{chowmotives}  
The comments on Thm.~\ref{poinrat} show we can express  the coefficients in the Poincar\'e series  from the trace of Frobenius iterates acting on  the $p$-adic cohomology that underlies Dwork's zeta rationality result. 
Positive: The computation is effective. Negative: The cohomology underlying Dwork's construction varies with $p$. Nothing in 0 characteristic represents it. 

Even, however, Dwork's cohomology \cite{Dw60} deals with stratifying your original variety. By \lq\lq combining\rq\rq\ the different pieces you conclude the rationality of the zeta function from information on the Frobenius action from the hypersurface case. The explicit estimates that Thm.~\ref{poinrat} relies on for these hypersurface computations are from \cite{Bo78}. 

Denef and Loeser \cite{DL98} applied Galois stratification (see the arXiv version of \cite[App.]{Ha07}) to eliminate quantifiers in their $p$-adic  problem goals. Their technique, as in \cite{DL02}, applies to consider -- for almost all $p$ -- how to express those Poincar\'e series, as stated in the abstract to this paper, as elements in the Grothendieck group generated by completely reducible (by the action of the Frobenius) pieces of the weighted $\ell$-adic cohomology -- twisted by Tate Modules; a tensoring of the group by some power of the cyclotomic character   -- of {\sl nonsingular projective\/} varieties.  Thus, the coefficients -- as in the previous  cases -- derive from the trace applied to restricting powers of the Frobenius for $p$ to these coefficients. \S\ref{del} has further clarification. 

They use Galois stratification, with the field crossing argument (Rem.~\ref{fieldcrossing}) doing a lot of work in putting covers and conjugacy classes to the background. Yet, this artistic enhancement  to the uniformity in $p$ in the uniform stratification (Comments on the proof of Thm.~\ref{poinrat}) is still akin to the previous methods.

\subsection{Deligne and motives} \label{del} The word   {\sl motive\/} refers to the weighted pieces -- rather than (pure) $m$th cohomology of a projective nonsingular variety -- being a summand of this, tensored by a Tate twist. A correspondence -- cohomologically idempotent --  is attached to indicate the source of the projection that detaches a summand from the pure weighted cohomology. Such idempotents are compatible with their appearance in writing the Poincar\'e series attached to a conjugacy class count in terms of characters induced from cyclic subgroups (say, \cite[\S 31.3]{FrJ04}$_1$; as in Thm.~\ref{poinrat}) of the group $G$. 

For one, from the main theorem of \cite{De74}, the absolute values of the eigenvalues of the Frobenius on these pieces are known. That allows more carefully considering any cancelation of these eigenvalues.  For example, \cite[Thm.~8.1]{De74} gives a definitive result on the eigenvalues of a complete intersection of dimension $n$. The error term there is $O(q^{n/2})$ by contrast to the usual expection that an error term is $O(q^{n\nm 1/2})$ as in \eqref{Frobest}. 

The \S\ref{schurcover} example has a naturally attached projective nonsingular variety, $X_2^\dagger$ to express the Poincar\'e series coefficients for a Galois stratification $\sS_0$ for (almost) all $p$, not just for those in the exceptional set. Even this example shows an aspect of using Denef-L\"oser that the previous two approaches don't seem to have: 
$$ \text{Using (Chow) motive pieces to relate uniform stratification primes and incidental primes.} $$ 

I've wanted to say this in print for a long time. It is common to think of $\ell$-cohomology as if it could all be from the cohomology of abelian varieties. For example, \cite{De72a} expresses the Weil conjectures for the cohomology of K3 surfaces (that in the middle being the only significant piece)  {\sl is\/}  through a Clifford algebra from the cohomology of an (nonobvious) abelian variety with appropriate Tate twists. Indeed, actual descent to give Frobenius action in positive characteristic requires a rigidying argument. This  uses endomorphisms from the Clifford algebra. So, the eigenvalue weights come from Weil's Theorem on abelian varieties. 

Yet -- has anyone considered the question in print earlier -- \cite{De72b} shows the middle cohomology of most complete intersections is not expressible from abelian varieties. 

\begin{quest} How, in general, would you distinguish those Chow motive pieces that do come from abelian varieties? \end{quest} 

The archetype for such a consideration is the classically considered use of minimal models to describe appropriate Hasse-Weil zeta functions of elliptic curves, one that includes factors for the primes of bad reduction. I say this even though Denef-L\"oser stratification replacement uses resolution of singularities in 0 characteristic. That is, their method won't directly touch the primes of the  incidental stratification. 

For good reason, I should consider how my viewpoint can properly and practically tackle the motivic aspects; also what  \cite{Ni07} gains using  Voevodsky's motives. 

\section{Motivations for extending Galois Stratification} \label{motivations}  

In each of the remaining subsections I put a classical veneer -- appropriate to areas of specific researchers  -- on papers that extend the Galois stratification procedure. 

\subsection{The author's motivations} \label{author} Very special problems motivated surprisingly general approaches to diophantine equations. Well-known examples drove the author's motivations and his many conversations with the principles involved. We've already mentioned Artin's Conjecture. One short breathless paragraph summarizes another. 

The intense work on zeta functions at {\sl all\/} primes for an elliptic curve, $E$, over $\bQ$ aimed at getting a functional equation for the associated Hasse-Weil zeta function $Z(E)$ that -- combined with the Eichler-Shimura congruence formula for the Frobenius on modular curves and the Shimura-Taniyama-Weil conjecture -- motivated the use of $Z(E)$ to conjecture -- a la Birch-Swinnerton-Dyer -- what is the rank of $E(\bQ)$ points. 

Indeed, less breathlessly, what natural diophantine statements could be expected to have such Euler factors at {\sl every\/} prime? Before Ax and Kochen, most number theorists seemed convinced of the naturalness of the Artin Conjecture: such a close analog of Chevalley's Theorem! It, however, no longer appears to be so natural. The Langlands Program is quite dependent on specifics -- of which there are few examples --  in the paragraph above. Yet the belief in Hasse-Weil zeta functions rest on the naturalness of the Langlands Program. 

\S\ref{questpoin} concludes -- via a cohomology component -- with a strong uniformity  statement, giving an object over $\bQ$ from which Euler factors for almost all primes $p$ come from reductions moduli those primes. Then, \S\ref{unifp}, based on Denef-L\"oser, strengthens the cohomology analogy with the work on elliptic curves. 

Finding semi-classical examples of such diophantine problems was this author's motivation. Exceptional covers, as in Ex.~\ref{coverpaircont}, is one of a general type, which like the Artin-Conjecture, starts as a  simple quantified statement. Generalizing the property of exceptional covers and {\sl Davenport pairs\/}  (as in \cite{Fr05} and \cite{Fr12}) tied many seemingly disparate diophantine considerations together. 

Toward the virtues of translating diophantine statements to Chow motives one property seems aesthetically significant: Separating zeta coefficients that are rational functions in Tate motives -- and therefore, on evaluating $\Fr_q$ at them, give rational functions in $q$  -- from those that are coefficients involving more complicated cohomology. 

In Prop.~\ref{exceptprop}, \eqref{exccondbp} expresses that for infinitely many Poincar\'e factors, the coefficients that appear are rational functions in Tate Motives. Expression \eqref{exccondbpnot} suggests that usually there will also be infinitely many Chow motive coefficients that are linear combinations of cohomology not composed from Tate motives.  

Rem.~\eqref{langlands}, about elliptic curves over $\bQ$ with the $\GL_2$ property in Serre's Open Image, poses one of the most simply stated unsolved problems for the Langlands program.  
This concludes my summary of my  many-years-ago motivations.  
 
 \subsection{Generalizing Felgner to Frobenius Vector problems}   \label{hrushwan} 
 
 Denote an algebraic closure of $\bF_q$ by $\bar \bF_q$.  We regard such a generalization as considering whether there could be analogous results in quantifying our variables according to Def.~\ref{frobvec}.  
\begin{defn} \label{frobvec} For a given prime-power $q$, and $\row d m$ an $m$-tuple of integers, refer to $\Fr_q^{\bd}\eqdef (\Fr_q^{d_1},\dots, \Fr_q^{d_m})$ as a {\sl Frobenius vector}. Then, $\Fr_q^{\bd}$ acts on $\bar \bF_q^m$ coordinate-wise, allowing us to speak of the elements $\bx\in \bar \bF_q^m$ fixed by a Frobenius vector. Denote these $\bar \bF_{q}^{m}(\Fr_q^\bd)$.  \end{defn} 

Notice we may write the elements -- in $\bF_{p^2}$, referenced by Felgner's problem, as $\bar \bF_{q}^{m}(\Fr_q^\bd)$ with $\bd=(2,\cdots,2)$. Frobenius vectors allow generalizing elementary statements (and Galois stratifications)  to where variables have values in differing extensions of $\bF_q$. With the notation above,  consider two Frobenius vectors: $\Fr_q^\bd$ of length $m$; and $\Fr_q^{\be}$ given by $\be=(\row \be k)$  with $\be_i$ of length $n_i$. Thus, $\be$  has  length $N_k=\sum_{i=1}^k n_i$. 

Here are some basic questions. 
\begin{edesc} \label{frobv} \item \label{frobva} Are the corresponding Poincar\'e series of Galois stratifications, with no quantified variables, rational? 
\item \label{frobvb}  Are there Bombieri-Dwork bounding degrees of the involved rational functions? 
\item \label{frobvc}  Does the Galois stratification procedure generalize to give rational functions?  That is, given quantiers, can we eliminate them to be at  \eql{frobv}{frobva}  with bounds from \eql{frobv}{frobvb}? 
\item \label{frobvd}  If the above, when are these Poincar\'e series new; not expressed from series associated to our previous Galois stratifications? 
\end{edesc}

\display{Wan's zeta functions} 
There are no Galois stratifications or quantified variables in \cite{W03}. Just the definition of a zeta function defined by a Frobenius vector. That is coefficients defined from counting points $\bx\in \bar \bF_q^m(\Fr_q^\bd)$ as above, a la Dwork, on an affine variety in $\afA^m$. 

Here are its contributions to \eql{frobv}{frobva}. Following a preliminary result by Faltings, documented in \cite{Wa01}, \cite[Thm.~1.4]{W03} shows the zeta  is  a rational function. 

Then it considers \eql{frobv}{frobvb} on the total degree of the zeta function akin to what Thm.~\ref{poinrat} uses. There is a preliminary result for special Frobenius vectors with $d_1|d_2| \cdots |d_m$ (consecutive $d\,$s dividing the next in line) in \cite{FuW03}  based on Katz's explicit bound for $\ell$-adic Betti numbers in \cite{Ka01}. Then,   \cite[Conj.~1.5]{W03} has a conjecture that there is an explicit total degree bound in general.  
 
\subsection{Logicians Chatzidakis, Hrushovski and Tomasovic}  \cite{CH99}, \cite{Hr12} and \cite{To16} deal with the theory of difference fields (or schemes), and the generalization of Galois formulas (akin to \cite{FrJ04}) over difference rings. Difference fields include the field $\bar \bF_q$ together with a symbol for the Frobenius automorphism. The generalization includes the addition of a symbol $\sigma$ for a distinguished automorphism of a scheme.

\display{Galois stratification abstracted} Recall, the main device of \cite{FrS76} was -- necessarily -- going beyond the original set of questions, to the richer set of {\sl Galois stratifications}. 
A short statement of what then happened: a quantifier elimination (as in \S\ref{nonregcheb}) gave a primitive recursive procedure for the theory of {\sl first-order definable sets in the language of schemes over finite fields\/} equipped with powers of the Frobenius automorphism (and related theories). 

That statement is similar to the \cite[p.~2, Thm.~1.1]{To16} main result. Except, the word {\sl twisted\/} appears in front of {\sl Galois Stratification\/}, and {\sl difference schemes\/} replace {\sl schemes\/}. 

Further, \cite[Thms. Thms. 1.1 and 1.2]{To16} are -- with these enhancements  -- the exact analog of results in  \cite{FrS76}, and proceed on the same basis. That is, you start with a Galois stratification over a {\sl difference\/} scheme $(X,\sigma)$. Here the data is a stratification of $(X,\sigma)$ into difference subscheme pieces $(X_i,\sigma), i\in I$, where each $X_i$ is equipped with an \'etale Galois cover $Z_i\to X_i$ with group $G_i$ and conjugacy classes $\bfC_i$ in $G_i$. 

Defining, though, appropriate Galois covers has subtleties about extending $\sigma$ and assuring with the additional data that \lq\lq difference field specializations\rq\rq\  of points of $X_i$ lift to corresponding extending specializations in the Galois cover. 

The idea is just like a Galois formula attached to affine space based on achieving Frobeniuses in the conjugacy classes attached to points of $X$ running over finite fields, subject to quantifiers placed on some of the variables, or for counting such achievements. 

\display{Additional complications} More general here is that $\sigma$ is added structure on $X$ from an endomorphism of $X$. Also,  everything is in the language of difference schemes. Many notions aren't defined in this paper, but in a longer paper called \lq\lq Twisted Lang-Weil\rq\rq\ \cite{To14}. 

The essence of \cite{To16} is to define the direct image, $\sB$, of the twisted Galois Stratification $\sA$, given an \'etale morphism $(X,\sigma) \to (Y,\sigma')$, so that their underlying Galois formulas are appropriately related: an exact analog of \cite{FrS76} and \S \ref{quoteFrJ04}.  

Elimination of quantifiers follows from applying this to projection of affine $n$-space onto affine $m$-space $m< n$. Until, inductively, you are at a quantifier free statement on a (almost certainly) complicated Galois stratification. Complicated or not,  a Chebotarev analog applies to decide if almost all Frobeniuses end up in the associated conjugacy classes. 

The difficulty is in extending the definitions to include difference schemes in both pieces. Defining a Galois cover of difference schemes \cite{To14}, so that there is a Chebotarev analog, and then proving the endomorphisms  survive, in appropriate shape, the process of direct image. 

Both apply to \eql{frobv}{frobvc} on eliminations of quantifiers, though neither goes after the Poincar\'e series. There is even the notion of an existentially closed difference field, \cite{CH99}, the analog of the (existentially closed -- \PAC) Frobenius fields of \S\ref{quoteFrJ04} as in \eql{frfi}{frfia}.

Certainly, the very long \cite{Hr12} preprint and \cite{To16} cover territory akin and in ways more general than this paper, though without the attempt to consider the relation with zeta functions. To give a hint as to what it does include, we conclude with a few comments on it starting with the {\sl Twisted Lang-Weil Theorem}. Look back at \S\ref{nonregcheb} for the extensive use made of variants on the original Lang-Weil Theorem. 

Here is an example of what Hrushovsky and Toma\v si\c c are considering.  

\begin{exmp} Suppose $\phi: A\to A$ is a nontrivial endomorphism of an abelian variety over $\bF_q$. Denote its graph by $S$. 

With $\Gamma_{\Fr_q}$ the graph of the Frobenius map $\ba \mapsto \ba^q$, consider $|S(\bar \bF_q)\cap \Gamma_{\Fr_q}(\bar \bF_q)|$. The cardinality of this set is the same as that of the elements in $\bar \bF_q$ for which $\phi(\ba)=\ba^q$. \end{exmp} 

In generality, replace $A$ by any variety $X$ and replace $S$ by any correspondence with finite projections to each factor of $X\times X$ (there are details to state this completely correctly). The nonobvious estimate for this cardinality is essentially the same as for Lang-Weil, and both authors recount the history behind the result initiated (apparently) by a conjecture of Deligne in \cite{De74}, with an ultimate proof dependent on his proof of the Weil conjectures. The case with  $S$ the diagonal exactly gives Lang-Weil.  

As \cite{Hr12} exposits, this is a key ingredient in the theory of $\bar \bF_q$ with its Frobenius. Again, they are each thinking of the generalization to difference schemes for which the key word is endomorphism of $X$. So, \cite{To16} is considering "Galois covers" $\hat \phi: \hat C \to X$ with $X$ equipped with an endomorphism $\sigma$ and $\hat C$ equipped with endomorphisms $\Sigma$, above it, that can be regarded as closed under something akin to conjugation. For this situation they must prove the analog of the Chebotarev Theorem counting the number of times the Frobenius is achieved within $\Sigma$ for points $X(\bar \bF)$. 

For logicians the key question: Does this gives a theory of various kinds of difference fields. That would include $\bF_q$ with the action of the Frobenius $\Fr_q$? The answer is Yes! That is what their papers show. In the Frobenius case, how does that work with our theme of uniformity in $p$,  an algebraic-geometric property sort-of compatible  with their logic frameworks? 

\begin{quest} What has this to do with Denef-L\"oser and with Wan? \end{quest} 

\!\!\!\!Is this difference equation approach  like the production of Poincar\'e series with those Chow motive coefficients?  Especially since Felgner's simple question offered so many peeks at the value of a decision procedure. 

For example, \cite{Hr12} considers specific diophantine applications to a very interesting field: The cyclotomic closure, $\bQ_\cyc\eqdef \cup_{n=1}^\infty \bQ(e^{\frac{2\pi i}n})$. 
This field is especially interesting because it is a characteristic 0 analog of the algebraic closure of a finite field. Consider its mysteries in light of the Fried-V\"olklein conjecture \cite{FrV92}, a strengthening of its main Theorem. 

\begin{guess} If  $K\le \bar \bQ$ is Hilbertian (Hilbert's irreducibility theorem holds in $K$) and its absolute Galois group, $G_K$,  is projective (among profinite groups), then $G_K$ is pro-free.\end{guess}  

The special case $\bQ_\cyc$ is a conjecture of Shafarevich. If true, in lieu of the theory of Frobenius fields, this connects $\bQ_\cyc$ in an even deeper way with the theory of finite fields. That, by the way, is a big piece of the motivation for \cite{FrJ04}. 

\subsection{The author states that there is no conflict of interest} 
 
 \providecommand{\bysame}{\leavevmode\hbox to3em{\hrulefill}\thinspace}
\providecommand{\MR}{\relax\ifhmode\unskip\space\fi MR}
\providecommand{\MRhref}[2]{%
\href{http://www.ams.org/mathscinet-getitem?mr=#1}{#2}}
\providecommand{\href}[2]{#2}
\!\!\!

\end{document}